\documentclass[9pt]{article}
\usepackage{epsfig}
\usepackage{amsmath,amsthm,amsfonts,amssymb,amscd,euscript}
\usepackage[numbers,sort&compress]{natbib}
\usepackage{amsmath}
\usepackage{amsthm}
\usepackage{amssymb}
\usepackage{latexsym}
\usepackage[all]{xy}
\usepackage{tikz}

\theoremstyle{definition}
\newtheorem{definition}{Definition}[section]

\newtheorem{example}[definition]{Example}

\newtheorem{remark}[definition]{Remark}

\theoremstyle{plain}
\newtheorem{theorem}[definition]{Theorem}
\newtheorem{lemma}[definition]{Lemma}
\newtheorem{proposition}[definition]{Proposition}
\newtheorem{corollary}[definition]{Corollary}

\begin{document}
\title{On Semi-simplicity Results in Residuated Lattices}
\author{Esmaeil  Rostami$^1, ^2$\\ \footnotesize
	 $^1$ Mahani Math Center, Afzalipour Research Institute, Shahid Bahonar University of Kerman, Kerman, Iran \\ \footnotesize  $^2$ Department of Pure Mathematics, Faculty
	 of Mathematics and Computer, Shahid
	 Bahonar University of Kerman, Kerman,
	 Iran\\ \footnotesize Email: e\underline{ }rostami@uk.ac.ir}

\maketitle

\begin{abstract} We develop the theory of residuated lattices by introducing and studying several new types of filters and related concepts, including semi-simple filters, essential filters, the socle of a filter, and independent families of filters.  Our primary goal is to understand the inner structure of residuated lattices by analyzing these new objects.
First, we establish the key properties of simple and essential filters. Next, we then provide both algebraic and topological characterizations for identifying when a filter is simple or essential.
Furthermore, we use the concepts of the socle and independent families to consider deeper into the structure of filters and the residuated lattice itself. We also provide  several characterizations for semi-simple filters and residuated lattices. A central result shows that for finite residuated lattices, being semi-simple is equivalent to being hyperarchimedean, highlighting the natural connection between these concepts. Complementary results deepening on the understanding of the relation between simple filters and maximal filters in residuated
lattices are also established.\end{abstract}


\markboth {Esmaeil Rostami}
{Pseudo-irreducibility and Boolean lifting property}

\section*{Introduction}
Algebraic structures play a fundamental role in formalizing different logical systems. Boolean algebras, for instance, are the natural framework for classical logic. However, to study non-classical logics such as fuzzy logic where the logical \lq\lq and" can depend on the order of propositions we need more flexible tools. Residuated lattices serve as this primary, general framework.

Introduced by Ward and Dilworth in 1939 [17], residuated lattices generalize the ideal theory of rings. Important subclasses like MV-algebras and BL-algebras have a direct and well-established connection to fuzzy logic [5, 7].

In parallel, the concepts of simplicity and semi-simplicity are central in ring and module theory. A simple ring has no non-trivial two-sided ideals, while a semi-simple ring is, by the fundamental Wedderburn-Artin Theorem, a finite direct product of simple rings. This decomposition provides a powerful structural insight [2,11]. The theory of semi-simple structures has since been extended beyond rings, for instance, the basic structure of semi-simple modules is detailed in [11], and the concept for semi-rings, based on ideals and coideals, has been introduced and investigated for its homological properties [3, 8].

However, the classical ring-theoretic perspective can sometimes obscure the underlying order-theoretic relationships. We argue that the language of residuated lattices  theory can offer a clearer and more unified viewpoint. In this lens, a simple module corresponds to a two-element chain, the simplest non-trivial lattice. A semi-simple module then corresponds to a Boolean algebra, which is a direct product of these simple chains. This approach provides an intrinsic, order based description of semi-simplicity and allows us to leverage the powerful tools of residuated lattice theory to gain a deeper understanding of these algebraic systems. This paper aims to develop this perspective systematically.

The paper is organized as follows. In Section \ref{1c}, we define simple and essential filters and investigate their key properties. Sections \ref{2c} and \ref{3c} are devoted to introducing independent families of filters and the socle of a filter, complete with characterizations using the maximal spectrum. In Section \ref{4c}, we define semi-simple filters. A central result we prove is that in finite residuated lattices, being semi-simple is equivalent to being hyperarchimedean, underscoring the fundamental nature of semi-simple structures.

\section{Preliminaries}  
In this section, we review some definitions and results which will be used throughout this paper.

\begin{definition}\label{2.111}[5, 7, 13] A \textit{residuated lattice}  is an algebra $(L, \wedge, \vee, \odot, \rightarrow, 0, 1)$ 
	of type $(2, 2, 2, 2, 0, 0)$ satisfying the following axioms:\\
	(RL1) $(L, \wedge, \vee, 0, 1)$ is a bounded lattice (whose partial order is denoted by $\leq $);\\
	(RL2) $(L, \odot, 1)$ is a commutative monoid;\\
	(RL3) For every $x, y, z \in L$, $x \odot z \leq y$ if and only if $z\leq x\rightarrow y$ (residuation).
\end{definition}
For   $x, y \in L$  and  $n\in \mathbb{N}$, we define: \[x^{\ast}:= x\rightarrow 0,  x^{**}:=(x^*)^*, x^0:=1,  \text{  and }  x^n:=x^{n-1}\odot x.\]

An element $x$  of a residuated lattice $L$ 
is called \textit{complemented} if there is an element $y \in L $ such that $x\wedge y = 0$ and $x\vee y= 1$, if such an element $y$ exists it is unique  and it is called  \textit{the  complement of} $x$. We  denote the  complement  of $x$ by $x'$.
The set of all complemented elements in $L$ is denoted by $B(L)$ and
is called \textit{the Boolean center of}
$L$.

In the following proposition, we collect some main and well-known properties of residuated lattices, and throughout the paper, we frequently use them without referring.

\begin{proposition}\label{2.222} [5, 7, 13] 
	Let $L$ be a residuated lattice, $x, y, z\in L$ and $e, f\in B(L)$. Then we have the following statements:\\
	\begin{enumerate}
		
		\item $x \leq y$ if and only if $ x \to y = 1$;
		\item   If $x \leq y$, then $y^*\leq x^*$; 
		\item  $x\odot y=0$ if and only if $x\leq y^*$. Also, $x\odot x^* = 0$;
		\item $x\odot(x\to y)\leq y$ and $y\leq x\to y$;
		
	\item If $x\vee y=1$, then $x\wedge y=x\odot y$ and $x^n\vee y^m=1$ for each $m, n\in \mathbb{N}$;
		\item $x\vee(y\odot z)\geq (x\vee y)\odot (x\vee z)$, $(x\vee y)^* = x^* \wedge y^*$, and $x\odot(y\vee  z)= (x\odot y)\vee (x\odot z)$;
	
		\item $e' = e^*$,  $e\wedge e^*=e\odot e^*=0$,  $e^{**} = e$,  $e\odot f= e \wedge f\in B(L)$, and for each $n\in \mathbb{N}$, $e^n=e$;
	
		\item $e\odot x = e \wedge x$, $e \wedge (x \odot y) = (e \wedge x) \odot (e \wedge y)$, and
	 $e \vee (x \odot y) = (e \vee x) \odot (e \vee y)$;

		\item $e\to x=e^*\vee x$, $x\to e=x^*\vee e$, $e\odot(e\to x)=e\wedge x$ and  $x\odot(x\to e)=e\wedge x$;

		\item If $x\vee x^*=1$, then $x\in B(L)$;
	\item $e\leq x$ if and only if $e\to x=1$ if and only if $e^*\vee x=1$.
	\end{enumerate}
	
\end{proposition}


\begin{definition}
	A non-empty subset $F$  of a residuated lattice $L$ is called a \textit{filter} 
	of $L$ if the following conditions hold:\\	
	(F1) If $x, y\in F$, then $x\odot y \in F$;\\
	(F2) If $x\leq y$ and $x\in F$, then $y \in F$.
\end{definition}

\begin{definition}
	A non-empty subset $F$  of a residuated lattice $L$ is called a \textit{deductive system} 
	of $L$ if the following conditions hold:\\	
	(D1) $1 \in F$;\\
	(D2) If $x, x\to y\in F$, then $y \in F$.
\end{definition}
A non-empty subset $F$ of a residuated lattice $L$ is a filter if and only if it is a deductive system, see [13] for more details.

We denote by $ Filt(L)$ the set of all filters of $L$.
A filter $F$ is called \textit{proper} if $F \neq L$. Clearly, a filter $F$ of $L$ is proper if and only if $0\not\in F$, equivalently, if for each  $x\in L$ we have either $x\not\in F$ or $x^*\not\in F$.

For a non-empty subset $S$ of a residuated lattice  $L$, we set $[S):=\bigcap\{F\in  Filt(L)\mid S\subseteq F\}$,   called \textit{the filter of $L$ generated by $S$}. We denote by $[x)$ \textit{the filter of $L$ generated by} $\{x\}$. A filter $F$ of a residuated lattice $L$ is called a \textit{finitely generated} filter if $F =[x_1, ..., x_n)$ for some $x_1, ..., x_n\in F$ and $1\leq n$.

 It is well-known that the  lattice $( Filt(L),\subseteq)$ is  complete. Actually,  for a family $\{F_i\}_{i\in \Lambda}$ of filters of $L$ we have $\bigwedge_{i\in \Lambda} F_i=\bigcap_{i\in \Lambda}F_i $ and $\bigvee_{i\in \Lambda} F_i=[\bigcup_{i\in \Lambda}F_i )$. If $F\in  Filt(L)$, we set $F^*:=\{x\in L\mid F\cap [x)= \{1\}\}$, see [13, 15].

\begin{proposition}\label{2.555}\label{3.7} [13] Let  $S$ be a non-empty subset of a residuated lattice $L$, $x, y \in L$, $F, G\in  Filt(L)$, and $\{F_i\}_{i\in\Lambda}\subseteq  Filt(L)$. Then
	
	\begin{enumerate}
		\item $[S) = \{x\in L \mid  s_{1} \odot \cdots \odot s_{n}\leq x $ \text{  }\text{for some}\text{  } $n \geqslant 1 \text{  }\text{and}\text{  } s_{1},...,s_{n}\in S\}$;
		\item $[x) = \{z\in L \mid x^n\leq z $ for some $n \geqslant 1\}$. In particular, if $e, f\in B(L)$, then $[e) = \{z\in L \mid e\leq z \}$, and $[e) =[f)$ if and only if $e=f$;
		\item $\bigvee_{i\in \Lambda} F_i=[\bigcup_{i\in \Lambda}F_i )= \{x\in L \mid  $\text{  }\text{there are}\text{  }$ s_{1}\in F_{i_1},...,s_{n}\in F_{i_n} $\text{  }\text{such that}\text{  }$ s_{1} \odot \cdots \odot s_{n}\leq x $ \text{  }\text{for some}\text{  } $n \geqslant 1\text{  }\text{and}\text{  } s_i\in \Lambda\}$;
		\item  $[x)\vee [y)= [x \wedge y) = [x \odot y)$, and $[x)\cap [y) = [x \vee y)$;
	
		\item $F \vee G= [F \cup G) =\{x \in L \mid  a\odot b\leq x$ for some  $a\in  F$ and $b\in G\}$;
		\item $F^*\in Filt(L)$ and $F^*$ is the maximum filter of $L$ with respect to the property that $F\cap F^*=\{1\}$, that is, for $G\in Filt(L)$, $F\cap G=\{1\}$ if and only if $G\subseteq F^*$;
	\item $F^*=\{x\in L\mid x\vee y=1 \text{  }\text{for all}\text{ } y\in F\}$;
	\item The lattice $( Filt(L),\subseteq)$ is distributive.
	\end{enumerate}
\end{proposition}



Recall  that a proper filter $P \in  Filt(L)$ is called \textit{prime}
if  for $x, y\in L$, $x\vee y\in P$ implies either $x\in P$ or $y \in P$. We denote by $  Spec(L)$ the set of all prime filters of $L$. An easy argument shows that a proper filter $P$ is prime if and only if  for $F, G\in  Filt(L)$ if $F\cap G\subseteq P$ then we have either $F\subseteq P$ or $G\subseteq P$. 	Also, a proper filter $M\in  Filt(L)$ is called \textit{maximal} if $M$ is not strictly contained in a proper filter
of $L$. We denote by $Max(L)$ the set of all maximal filters of $L$. Clearly, every maximal filter  of  a 
residuated lattice is prime. Also, every proper filter is contained in a maximal filter, see [13, 15] for more details. 

\begin{theorem}\label{tpd}[13]
Let $M$ be a proper filter of a residuated lattice $L$. Then the following are equivalent:
\begin{enumerate}
	\item $M$ is a maximal filter of $L$;
	\item For any $x\in L$, $x\not\in M$ if and only if $(x^n)^*\in M$ for some $n\in\mathbb{N}$.
\end{enumerate}
\end{theorem}

Recall that a residuated lattice $L$ is called \textit{local} if $L$ has a unique maximal filter. Also, a residuated lattice $L$ is called \textit{semi-local} if $L$ has only a finite number of maximal filters.

The intersection of all maximal filters of a residuated lattice $L$ is called \textit{the radical
of L}, and will be denoted by $ Rad(L)$.

\begin{theorem}\label{rad}[13]
	Let $L$ be a residuated lattice. Then \[ Rad(L)=\{x\in L\mid  \text{for any}\text{  } n\in\mathbb{N}\text{  }\text{there is}\text{  } k_n  \text{  }\text{such that} \text{  }((x^n)^*)^{k_n} = 0\}.\]
\end{theorem}

Recall that  a  filter $P$ is called \textit{minimal prime} if $P \in   Spec(L)$ and, whenever $Q \in   Spec(L)$ and $Q \subseteq P$, then  we have $P = Q$. We denote by $ MinP(L)$ the set of all minimal prime filters of $L$.
\begin{theorem}\label{minp}[13]
Let $P$ be a prime filter of a residuated lattice $L$. Then the following are equivalent:
\begin{enumerate}
	\item $P$ is a minimal prime filter of $L$;
	\item For any $x \in P$ there is
	$y\in L\setminus P$ such that $x\vee y=1$.
\end{enumerate}	

\end{theorem}

Recall that a set $X$ with a family $\varGamma$ of its subsets
is called a \textit{topological space}, denoted by $(X, \varGamma)$,
if $\varnothing, X\in \varGamma$, the intersection of any finite number of members of $\varGamma$ is in $\varGamma$, and the arbitrary union of members of $\varGamma$ is in $\varGamma$. The members of
$\varGamma$ are called \textit{open} subsets of $X$, and the complement
of an open subset of $X$ is said to be  \textit{closed}.

If $Y$ is a subset of $X$, the smallest closed subset of $X$ containing $Y$ is called the \textit{closure} of $Y$ and denoted
by $\bar{Y}$, also the largest open subset of $X$ contained in $Y$ is called the \textit{interior} of $Y$ and denoted
by $Int(Y)$.

A point $x$ of a topological space $X$ is called an \textit{isolated point}, whenever the set $\{x\}$ is an open set in the topological space $X$.

For every subset $A$ of a residuated lattice $L$, we set $U (A) := \{P \in   Spec(L) \mid A \not\subseteq P\}$, and for each $a\in L$, we set $U (a):=U (\{a\})$.  The family $\{U(A)\}_{A\subseteq L }$  satisfies the axioms for open subsets for
a topology  over $  Spec(L)$. This topology is called  \textit{the Stone topology}. Since every maximal filter  is prime, we can consider $Max(L)$ as a subspace of $  Spec(L)$. For each $A\subseteq L $, we define $U_{Max}(A):=U(A)\cap Max(L)$. Then the  family $\{U_{Max} (A)\}_{A\subseteq L }$ satisfies the axioms for open subsets for
the subspace  topology over $Max(L)$.
For every subset $A$ of a residuated lattice $L$, we set $V (A) :=    Spec(L)\setminus U(A)$ and $V_{Max}(A):=Max(L)\setminus U_{Max}(A)$. Then the family $\{V(A)\}_{A\subseteq L }$  satisfies the axioms for closed subsets for
the Stone topology  over $  Spec(L)$, and   the  family $\{V_{Max} (A)\}_{A\subseteq L }$ satisfies the axioms for closed sets for the subspace topology over $Max(L)$, for more information see [6].

\begin{proposition}\label{2.700}[13] 
	Let $L$ be a  residuated lattice, $F\in  Filt(L)$ and $a\in L\setminus F$. Then we have the following statements:
	\begin{enumerate}
		\item There exists a prime filter $P$ of $L$ containing $F$ such that  $a\not\in P$;
		\item $F$ is the intersection of all prime filters which contain $F$, that is,  $F=\bigcap V(F)$;
		\item $\bigcap  Spec(L)=\bigcap MinP(L)=\{1\}$.
		
	\end{enumerate}
	
\end{proposition}

\begin{proposition}\label{smss} 
	Let $L$ be a residuated lattice, $F, G\in Filt(L)$ and $x, y\in L$. Then the following statements hold:
	\begin{enumerate}
		\item $F = L$ if and only if $V_{Max}(F)=\varnothing$ if and only if $V(F)=\varnothing$;
		\item $V(F) \subseteq V(G)$ if and only if $U(G) \subseteq U(F)$  if and only if $G\subseteq F$. In particular,  $V(F) = V(G)$ if and only if $U(F) = U(G)$ if and only if $F=G$;
	\item If $G\subseteq F$, then $V_{Max}(F)\subseteq V_{Max}(G)$ and $U_{Max}(G)\subseteq U_{Max}(F)$. In particular, if $R(L)=\{1\}$, then $V(F) = V(G)$ if and only if $V_{Max}(F)=V_{Max}(G)$ if and only if $F=G$;
		\item  $U(x) \cap U(y) =U(x\vee y)$ and  $U_{Max}(x) \cap U_{Max}(y) =U_{Max}(x\vee y)$;
		\item  $V(x) \cup V(y) =V(x\vee y)$ and  $V_{Max}(x) \cup V_{Max}(y) =V_{Max}(x\vee y)$;
	\item If $\{X_k\}_{k\in K}$ is a family of subsets of $L$, then $U(\bigcup_{k\in K}X_k)=\bigcup_{k\in K}U(X_k)$ and   $U_{Max}(\bigcup_{k\in K}X_k)=\bigcup_{k\in K}U_{Max}(X_k)$;
		\item If $\{X_k\}_{k\in K}$ is a family of subsets of $L$, then $V(\bigcup_{k\in K}X_k)=\bigcap_{k\in K}V(X_k)$ and   $V_{Max}(\bigcup_{k\in K}X_k)=\bigcap_{k\in K}V_{Max}(X_k)$.
	\end{enumerate}
\end{proposition}	
 An element $x$ of a residuated lattice $L$ is called
\textit{archimedean} if 
there is $n\in\mathbb{N}$  such that $x^n\in B(L)$. Also,
a residuated lattice $L$ is called
\textit{hyperarchimedean} if all its elements are archimedean.
\begin{theorem}\label{hyp}[13]
	For a residuated lattice $L$ the followings are equivalent:
	\begin{enumerate}
		\item $L$ is hyperarchimedean;
		\item $  Spec(L)=Max(L)$;
		\item $  Spec(L)= MinP(L)$.
	\end{enumerate}
\end{theorem}

\begin{theorem}\label{hyp1}[13] 
	If $L$ is  a hyperarchimedean residuated lattice, then $ Rad(L)=\{1\}$.
\end{theorem}
Recall that an element $a$  in a Boolean algebra is called a \textit{co-atom} if $a\not=1$
and if $a\leq b$, then we have either $a=b$ or $b=1$. Also, a Boolean algebra is said to be \textit{co-atom-less}, whenever it contains no co-atoms.
\section{Simple and Essential Filters}\label{1c}
We begin this section by recalling the definition of  a simple filter in a residuated lattice.

\begin{definition}
	A filter $F$ of a residuated lattice $L$ is called \textit{simple}, whenever  $F\not=\{1\}$ and if $G\subseteq F$ for some $G\in Filt(L)$, then we have either $G=\{1\}$ or $G=F$. We denote by $  SimpF(L)$  the set of all simple filters of a residuated lattice $L$.
\end{definition}
\begin{remark}\label{key}
	Let $F$ be a simple filter of a residuated lattice $L$. If $1\not=x\in F$, then since $\{1\}\not=[x)\subseteq F$, we have $F= [x)$. Hence every simple filter is principal and is generated by each $x\in F$ with $x\not=1$.
\end{remark}
\begin{remark}\label{t6}
	Let $L$ be a non-trivial residuated  lattice such that $ Filt(L)$ is a finite set (e.g., finite residuated lattices). Then clearly every non-trivial filter of $L$ contains a simple filter. But in general this fact is not true. For example assume that $L$ is a co-atom-less Boolean algebra (e.g., the Boolean algebra of clopen subsets of the Cantor Set, the Boolean algebra of Lindenbaum–Tarski Algebra of Propositional Logic or the Boolean algebra of measure algebra on $[0,1]$). Now, if $F$ is a non-trivial filter of $L$, then $F\not=\{1\}$. Hence there exists $1\not=a\in F$. Since $a$ is not a co-atom, there exists $b\in L$ such that $a
\lneq b \lneq 1$.

 Hence we have $\{1\}\not=[b)
\subsetneq [a)\subseteq F$ by Proposition \ref{2.555}(2). Thus $F$ is not a simple filter of $L$, and so $L$ has no simple filters. 
\end{remark}
In Theorem \ref{3.7}, we provide characterizations of residuated lattices $L$ such that every non-trivial filter of $L$ contains a simple filter.

According to the definition of a simple filter, it appears that this concept is dual to the notion of a maximal filter. 
In the following proposition, we consider the relation between simple and maximal filters that are generated by Boolean elements.
\begin{proposition}\label{smb}
	Let $e\in B(L)$ for a residuated lattice $L$. Then $[e)$ is a maximal filter of $L$ if and only if  $[e^*)$ is a simple filter of $L$.
\end{proposition}
\begin{proof} 
	$\Rightarrow)$ Let $[e)$ be a maximal filter of $L$. If $[e^*)=\{1\}$, then $e^*=1$, and so $e=0$. Thus $[e)=L$, which contradicts the maximality of $[e)$. Now if $\{1\}\subsetneq F\subseteq [e^*)$ for some filter $F$ of $L$. Choose $1\not=x\in F$. If $x\in [e)$, then $x\in[e)\cap [e^*)=[e\vee e^*)=\{1\}$, and so $x=1$, a contradiction. So $x\not\in [e)$ and hence $[x)\not\subseteq [e)$. Thus by Proposition \ref{2.555}(4) and the maximality of $[e)$, we have $[e\odot x)=[e)\vee[x)=L$. Hence by Proposition \ref{2.555}(1), there exists $n\in\mathbb{N}$ such that $e\odot x^n=e^n\odot x^n=(e\odot x)^n=0$. This shows that $x^n\leq e^*$, and so $e^*\in [x)$. Hence $[e^*)\subseteq [x)$. Now since $[x)\subseteq F$ and $F\subseteq [e^*)$, we have $F= [e^*)$. Therefore, $[e^*)$ is a simple filter of $L$.
	
	$\Leftarrow)$ Let $[e^*)$ be a simple filter of $L$. If $[e)=L$, then by Proposition \ref{2.555}(1) there exists $n\in\mathbb{N}$ such that  $e^n=0$, and so $e=0$. Thus $e^*=1$, and so $[e^*)=\{1\}$, which is impossible. Now if $[e)\subsetneq F$ for some filter $F$ of $L$. Choose $x\in F\setminus [e)$. Since $[x)\cap [e^*)\subseteq [e^*)$, we have either $[x)\cap [e^*)= [e^*)$ or $[x)\cap [e^*)=\{1\}$.
	
	If  $[x)\cap [e^*)= [e^*)$, then  $ [e^*)\subseteq [x)\subseteq F$. Thus $e, e^*\in F$. This shows that $0=e\odot e^*\in F$, and so $F=L$, in this case. If $[x)\cap [e^*)=\{1\}$, then $[x\vee e^*)=[x)\cap [e^*)=\{1\}$ by Proposition \ref{2.555}(5). This shows that $x\vee e^*=1$. Hence by Proposition \ref{2.222}, we have \[e=e\odot 1=e\odot (x\vee e^*)=(e\odot x)\vee(e\odot e^*)=e\odot x\leq x.\] Thus $e\leq x$, and so $x\in [e)$, a contradiction. Therefore, $[e)$ is a maximal filter of $L$.
\end{proof}

The following result establishes some properties of simple filters.
\begin{proposition}\label{133}
	Let $T$ be a simple filter of a residuated  lattice $L$. Then one of the following cases happens:
	\begin{enumerate}\label{stc}
		\item  There exists $x\in B (L)$ such that $T=[x)$. In this case, $[x^*)$ is the only maximal filter  of $L$ not containing $T$;
		\item For all $x\in T$, we have $[x^*)=L$. In this case, $T$ is contained in all maximal filters of $L$, that is, $T\subseteq  Rad(L)$.
	\end{enumerate}
\end{proposition}
\begin{proof} 
	Let  $1\not=x\in T$ be arbitrary. Since $T$ is simple and $\{1\}\not=[x)\subseteq F$, we have $T=[x)$. By Proposition \ref{2.555}(5), $[x\vee x^*)=[x)\cap[x^*)\subseteq T$. Hence one of the following cases happens:
	
	$\boldsymbol{Case\text{ }1:}$ If $[x\vee x^*)=\{1\}$ for some $1\not=x\in T$, then $x\vee x^*=1$, and so  $x\in B (L)$  by Proposition \ref{2.222}. Now since $T$ is a simple filter of $L$, by Proposition \ref{smb}, $M:=[x^*)$ is a maximal filter of $L$. It is obvious that $T\not\subseteq M $. Now, if $N$ is a maximal filter of $L$ with $M\not=N$, then since $x\vee x^*=1\in N$,  we have either  $x\in N$ or $x^*\in N$. If $x^*\in N$, then $M\subseteq N$ and since $M$ is maximal, we have $M=N$, a contradiction. Hence for all maximal filter $N$ of $L$ with $N\not=M$, we have $x\in N $, and so $T=[x)\subseteq N$. Therefore, in this case, $[x^*)$ is the only maximal filter of $L$ not containing  $T$.
	
	$\boldsymbol{Case\text{ } 2:}$ If $[x)\cap[x^*)=[x\vee x^*)=T=[x)$ for all $x\in T$, then $[x)\subseteq [x^*)$ for all $x\in T$, and so by Proposition \ref{2.555}(2) for each $x\in T$ there exists $n_x\in\mathbb{N}$ such that $(x^*)^{n_x}\leq x$. Hence $(x^*)^{n_x+1}=0$ for all $x\in T$, and so for all $1\not=x\in T$, we have $[x^*)=L$. If $a\in T$, then for all $n\in\mathbb{N}$  we have $a^n\in T$. Hence $[(a^n)^*)=L$, and so by Proposition \ref{2.555}(2) there exists $k_n\in\mathbb{N}$ such that $((a^n)^*)^{k_n}=0$. Thus by Proposition \ref{rad}, $a\in  Rad(L)$, that is,  $a$ is contained in all maximal filters of $L$. Hence, $T$ is contained in all maximal filters of $L$ in this case.
\end{proof}
As a direct consequence of Proposition \ref{stc}, we have the following corollaries.

\begin{corollary}\label{sr1}
	Let $T$ be a  non-trivial filter of a residuated lattice  $L$ with $ Rad(L)=\{1\}$. Then $T$ is simple if and only if there exists $e\in B (L)$ such that $T=[e)$ and $[e^*)$
	is a maximal filter of $L$.
\end{corollary}
\begin{proof}
	Since $ Rad(L)=\{1\}$ and $T\not=\{1\}$, we have $T\not\subseteq  Rad(L)$. Hence by Proposition \ref{stc}, if $T$ is a simple filter of $L$, then there exists $e\in B (L)$ such that $T=[e)$ and $[e^*)$ is the only maximal filter  of $L$ not containing $T$. The converse is true by Proposition  \ref{smb}.
\end{proof}

\begin{corollary}\label{scro1}
	Let $T$ be a simple filter of a residuated  lattice $L$ and $x\in T$. Then $[x^*)\not=L$ if and only if $x\in B(L)\setminus\{1\}$ and $T=[x)$.
\end{corollary}
\begin{proof}
	$\Rightarrow)$ By the proof of Proposition \ref{stc}, if $[x^*)\not=L$ for some $x\in T$, then $[x\vee x^*)=\{1\}$, $x\in B (L)$ and $T=[x)$. In this case, if $x=1$, then $x^*=0$ and so $[x^*)=L$, a contradiction. Hence $x\in B(L)\setminus\{1\}$ and $T=[x)$.
	
	$\Leftarrow) $ Let $x\in B(L)\setminus\{1\}$ and $T=[x)$. If $[x^*)=L$, then by Proposition \ref{2.555}, there exists $n\in\mathbb{N}$ such that $(x^*)^{n}=0$. Now since $x\in B(L)$, we have $x^*\in B(L)$, and so $x^*=(x^*)^{n}=0$ by Proposition \ref{2.555}. Thus, $x=x^{**}=1$, which is impossible. Therefore, $[x^*)\not=L$.
\end{proof}

\begin{example}\label{333.2}Let $L = \{0, a, b, c,  1\}$ denote the residuated lattice with Hasse diagram depicted in the following figure, and define the operations $\odot$ and $\rightarrow$ via the accompanying tables, see  [6, Example 4.4]:		
	\begin{center}  
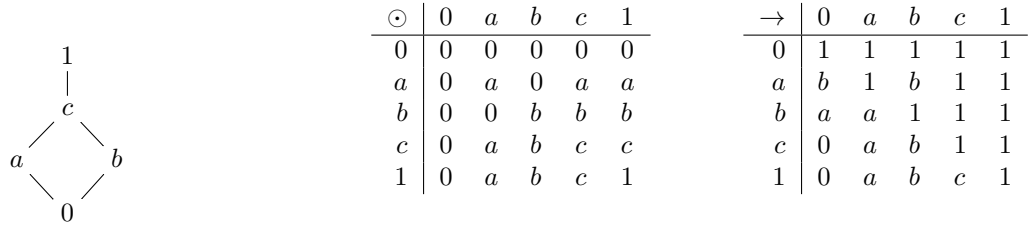
		\begin{minipage}{0.3\textwidth}
			\centering
			\[
			\xymatrix@=8pt{
				& & 1\ar@{-}[d]\\
				&   &c \ar@{-}[rd]\ar@{-}[ld] &    \\ 
				&a \ar@{-}[rd]& & b \ar@{-}[ld]\\
				&& 0}		\]
		\end{minipage}
		\begin{minipage}{0.65\textwidth}
			\centering
			\begin{tabular}{r|l l l l l }
				$\odot$ & 0 & $a$ & $b$ &$c$  &  1 \\\hline
				0& $0$ &$0$& $0$&$0$&$0$\\
				$a$ &$0$ &$a$ & $0$& $a$&$a$\\
				$b$ &$0$ &$0$& $b$&$b$&$b$\\
				$c$ & $0$ &$a$& $b$&$c$&$c$\\
				1 & $0$& $a$&$b$& $c$& $1$\\
			\end{tabular} 
			\hspace{1cm}
			\begin{tabular}{r|l l l l l }
				$\rightarrow$ & 0 & $a$ & $b$ &$c$  & 1 \\\hline
				0& $1$ &$1$& $1$&$1$&$1
				$\\
				$a$ &$b$ &$1$ & $b$& $1$&$1$\\
				$b$ &$a$ &$a$& $1$&$1$&$1$\\
				$c$ & $0$ &$a$& $b$&$1$&$1$\\
				1 & $0$& $a$&$b$& $c$& $1$\\
			\end{tabular}
			
			\vspace{0.5cm}
		\end{minipage}
	\end{center}
	\noindent	 The filter $[c)=\{c, 1\}$ is simple, but $c$ is not a Boolean element of $L$, and the filter $[c^*)=[0)=L$ is not maximal in $L$. Also, the filters $[a)$  and $[a^*)=[b)$ are maximal filters of $L$. Hence, Proposition  \ref{smb} is not true  in general, and the assumption $ Rad(L)=\{1\}$ cannot be omitted from Corollary \ref{sr1}.   Also, note that in this example, the filter $[c)$ is simple, but it is not prime or maximal filter, and the filter $[a)$  is maximal (and so it is prime), but it is not simple. Hence there are no general logical relationships between prime and simple filters.
\end{example}

\begin{definition}
	For a filter $F$ of a residuated lattice $L$, we set \[E_F:=\{H\in  Filt(L)\mid H\subseteq F \text{ and for each } G\in  Filt(L) \text{,  if }H\cap G=\{1\}  \text{,  then } F\cap G=\{1\}\}.\]

	When $H\in E_F$, we say that $H$ is an essential filter in $F$. In particular case, if $H\in E_L$, we say that $H$ is an essential filter in $L$.
\end{definition}

\begin{remark}\label{re}
	Clearly, for each filter $F$ of a residuated lattice $L$, we have $F\in E_F$. Hence $E_F\not=\varnothing$ for each filter $F$ of a residuated lattice $L$.
\end{remark}
\begin{remark}\label{eiclu}
	Let  $F, G, H\in  Filt(L)$ for  a residuated lattice $L$ such that $H\subseteq G\subseteq F$. Then, clearly, if $H\in E_F$, we have $H\in E_G$ and $G\in E_F$. Also, note that if $H\in E_G$ and $G\in E_F$, then $H\in E_F$.
\end{remark}

\begin{example}\label{555.2}Let $L = \{0, a, b, 1\}$ denote a lattice with Hasse diagram depicted in the following figure, and define the operations $\odot$ and $\rightarrow$ via the accompanying tables:		
	\begin{center}  
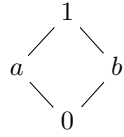
		\begin{minipage}{0.3\textwidth}
			\centering
			\[
			\xymatrix@=8pt{
				&   &1 \ar@{-}[rd]\ar@{-}[ld] &    \\ 
				&a \ar@{-}[rd]& & b \ar@{-}[ld]\\
				&& 0}		\]
		\end{minipage}
		\begin{minipage}{0.65\textwidth}
			\centering
			\begin{tabular}{r|l l l l }
				$\odot$ & 0 & $a$ & $b$   &  1 \\\hline
				0& $0$ &$0$& $0$&$0$\\
				$a$ &$0$ &$a$ & $0$& $a$\\
				$b$ &$0$ &$0$& $b$&$b$\\
				1 & $0$& $a$&$b$&  $1$\\
			\end{tabular} 
			\hspace{1cm}
			\begin{tabular}{r|l l l l  }
				$\rightarrow$ & 0 & $a$ & $b$   & 1 \\\hline
				0& $1$ &$1$& $1$&$1$\\
				$a$ &$b$ &$1$ & $b$& $1$\\
				$b$ &$a$ &$a$& $1$&$1$\\
				1 & $0$& $a$&$b$& $c$\\
			\end{tabular}
			
			\vspace{0.5cm}
		\end{minipage}
	\end{center}
	\noindent	 The filter $[a)=\{a, 1\}$ is not essential in $L$ since $[b)\cap [a)=\{1\}$ and $[b)\cap L\not=\{1\}$, but it is a  maximal (and so it is prime) and simple filter in $L$. Also, the filter $L$  is essential in $L$, but it is neither  prime (and it is not maximal) nor simple. Hence there is no general logical relationships between prime filters and essential filters. Also, there are no general logical relationships between simple filters and essential filters
\end{example}

\begin{example}\label{ddd}
	In Example \ref{333.2}, the filter $[c)$ is simple and   essential, but it is not prime.
\end{example}

\begin{example}\label{asa}
	In any non-trivial local residuated lattice, since every proper filter is contained in the unique maximal filter, we have the unique maximal filter is essential.
\end{example}

The following proposition is a direct consequence of the definition of $F^*$, Proposition \ref{2.555}, and Remark \ref{re}.
\begin{proposition}\label{esse}
	Let $F$ be a filter of a residuated lattice $L$. Then the following are equivalent:
	\begin{enumerate}
		\item $ F$ is an essential filter in $L$;
		\item $F^*=\{1\}$;
		\item For $x\in L$, if $x\vee y=1$ for all $y\in F$, then $x=1$;
		\item For each $1\not=x\in L$, there exists $y\in F$ such that  $x\vee y=1$.	
	\end{enumerate}
\end{proposition}

\begin{lemma}\label{efe}
	Let $F, H \in  Filt(L)$ for a residuated lattice $L$ such that $H\subseteq F$. Then $H\vee(H^*\cap F)\in E_F$.
\end{lemma}
\begin{proof}
	Since $H\subseteq F$, we have $H\vee(H^*\cap F)\subseteq F$. Now if $(H\vee(H^*\cap F))\cap G=\{1\}$ for some $G\in Filt(L)$,  then we have $(H\cap G)\vee((H^*\cap F)\cap G)=\{1\}$ by Proposition \ref{2.555}(8). Hence we have $H\cap G=\{1\}$  and $(H^*\cap F)\cap G=\{1\}$. From $H\cap G=\{1\}$, we conclude that $G\subseteq H^*$ by Proposition \ref{2.555}(6). Thus from $(H^*\cap F)\cap G=\{1\}$, we have $F\cap G=\{1\}$. Therefore, $H\vee(H^*\cap F)\in E_F$. 
\end{proof}
Recall that a subset $A$ of a topological space $X$ is called \textit{nowhere dense} if its closure has empty interior.
    
The following theorem states a topological characterization for a filter to be essential in a residuated lattice $L$ with $ Rad(L)=\{1\}$.
   
\begin{theorem}\label{111}
	Let $F$ be a filter of a residuated lattice $L$ with $ Rad(L)=\{1\}$. Then $F\in E_L$ if and only if $V_{Max}(F)$ is a nowhere dense subset of $Max(L)$.
\end{theorem}    

\begin{proof} 
	$\Rightarrow).$ Assume that $F\in E_L$. Since  $Int_{Max}(V_{Max}(F))$ is an open subset of $Max(L)$, there exists $G\in Filt(L)$ such that $Int_{Max}(V_{Max}(F))=U_{Max}(G)$. Assume that $U_{Max}(G)\not=\varnothing$. Choose $M\in U_{Max}(G)$. Hence $G\not\subseteq M$. Now if $\bigcap V_{Max}(G)\subseteq M$, then $G\subseteq M$, a contradiction. Hence, 	$\bigcap V_{Max}(G)\not\subseteq M$. Choose $x\in (\bigcap V_{Max}(G))\setminus M$, and hence $x\in N$ for all $N\in V_{Max}(G)$.
	
	Now assume that $y\in F$. Then $x\vee y\in N$ for all $N\in V_{Max}(G)$. If $N\in Max(L)\setminus V_{Max}(G)$, then $N\in U_{Max}(G)=Int_{Max}(V_{Max}(F))\subseteq V_{Max}(F)$, hence $F\subseteq N$, and so $y\in N$. Thus $x\vee y\in N$. Therefore, $x\vee y\in N$ for all $N\in Max(L)$, that is, $x\vee y\in Rad(L)=\{1\}$, and so $x\vee y=1$. From $x\not\in M$, we have $x\not=1$. Thus, $F\not\in E_F$ by Proposition \ref{esse}, which contradicts the assumption. Therefore, $Int_{Max}(V_{Max}(F))=U_{Max}(G)=\varnothing$, that is, the closed subset $V_{Max}(F)$ of $Max(L)$ is nowhere dense in $Max(L)$.
	
	$\Leftarrow).$ Assume that $V_{Max}(F)$ is a nowhere dense subset of $Max(L)$. Hence $Int_{Max}(V_{Max}(F))=\varnothing$. Let $1\not=x\in L$. Since $ Rad(L)=\{1\}$, we have $x\not\in  Rad(L)$. Hence, $U_{Max}(x)\not=\varnothing$. By assumption, $U_{Max}(x)\not\subseteq V_{Max}(F)$, and so $U_{Max}(x)\cap U_{Max}(F)\not=\varnothing$. Choose $M\in U_{Max}(x)\cap U_{Max}(F) $. Then $x\not\in M$ and $F\not\subseteq M$. Choose $y\in F\setminus M$. If $x\vee y=1$, then $x\vee y\in M$, and so we must have either $x\in M$ or $y\in M$ since every maximal filter is prime, which is impossible. Hence, $x\vee y\not=1$. Therefor for each $1\not=x\in L$, there exists $y\in F$ such that $x\vee y\not=1$. Hence $F\in E_F$ by Proposition \ref{esse}. 
\end{proof}
\begin{remark}\label{keyy}
	In Example \ref{333.2}, $V_{Max}([c))=Max(L)$ is not a nowhere dense subset of $Max(L)$, but $[c)\in E_L$. Hence the hypothesis $ Rad(L)=\{1\}$ in  Theorem \ref{111} is necessary. 
\end{remark}

\begin{proposition}\label{den}
	Let $L$ be a residuated lattice and $\mathfrak{A}\subseteq   Spec(L)$. Then the closure of $\mathfrak{A}$ in $  Spec(L)$ is equal to $V(\bigcap A)$. In particular, if $\mathfrak{A}\subseteq Max(L)$, then the closure of $\mathfrak{A}$ in $Max(L)$ is equal to $V(\bigcap A)\cap Max(L)$.
\end{proposition}
\begin{proof}
	Since $\mathfrak{A}\subseteq V(\bigcap \mathfrak{A})$ and $V(\bigcap A)$ is a closed subset of $  Spec(L)$, we have $\bar{\mathfrak{A}}\subseteq V(\bigcap \mathfrak{A})$. Now let $\mathfrak{C}$ be a closed subset of $  Spec(L)$ such that $\mathfrak{A}\subseteq \mathfrak{C}$. From the definition, there exists $F\in Filt(L)$ such that $\mathfrak{C}=V(F)$. Hence $\mathfrak{A}\subseteq V(F)$, and so for each $P\in \mathfrak{A}$, we have $F\subseteq P$. Thus $F\subseteq \bigcap \mathfrak{A}$, and so $V(\bigcap \mathfrak{A})\subseteq V(F)=\mathfrak{C}$ by Proposition \ref{smss}(2). This shows that $V(\bigcap \mathfrak{A})$ is the smallest closed subset of $  Spec(L)$ that contains $\mathfrak{A}$. Therefore, $\bar{\mathfrak{A}}=V(\bigcap \mathfrak{A})$. The last part of the proposition follows from the fact that if $ X$
	is a topological space and  $Y$
	is a subspace of $X$, then for
	a subset $A$ of $Y$ the closure of $A$ in $Y$ is equal to the intersection of $Y$ and the closure of $A$ in $X$, see [4, Proposition 2.1.1].
\end{proof}

\begin{corollary}\label{ddden}
	Let $L$ be a residuated lattice and $\mathfrak{A}\subseteq   Spec(L)$. Then  $\mathfrak{A}$ is a dense  subset of $  Spec(L)$ if and only if $\bigcap \mathfrak{A}=\{1\}$.
\end{corollary}
\begin{proof}
	By the definition and Propositions \ref{2.700}, \ref{smss}, and \ref{den}, we have 
	\begin{align*}
	\mathfrak{A} \text{ is dense in  }   Spec(L) \Leftrightarrow &  \text{  } \bar{\mathfrak{A}}=  Spec(L)\\ \Leftrightarrow & \text{  } V(\bigcap \mathfrak{A})=  Spec(L) \\ \Leftrightarrow &\text{  } V(\bigcap \mathfrak{A})=V(\{1\}) \\ \Leftrightarrow &  \text{  } \bigcap \mathfrak{A}=\{1\}.
	\end{align*}
\end{proof}
\begin{corollary}\label{dddenm}
	Let $L$ be a residuated lattice such that $ Rad(L)=\{1\}$ and $\mathfrak{B}\subseteq Max(L)$. Then  $\mathfrak{B}$ is a dense  subset of $Max(L)$ if and only if $\bigcap \mathfrak{B}=\{1\}$.
\end{corollary}
\begin{proof}
	Since $ Rad(L)=\{1\}$, for a filter $F$ of $L$ we have  $V_{Max}(F)=Max(L) $ if and only if $ V_{Max}(F)=V_{Max}(\{1\}) $ if and only if $F=\{1\}$, see Proposition \ref{smss}(2). Hence by Proposition \ref{den}, we have 
	\begin{align*}
	\mathfrak{B} \text{ is dense in  } Max(L) \Leftrightarrow &\text{  }
	\text{The closure of }\mathfrak{B} \text{  in  } Max(L) \text{ is equal to } Max(L) \\ \Leftrightarrow & \text{  } V(\bigcap \mathfrak{B})\cap Max(L)=Max(L)\\ \Leftrightarrow &  \text{  }V_{Max}(\bigcap \mathfrak{B})=Max(L) \\ \Leftrightarrow & \text{  }V_{Max}(\bigcap \mathfrak{B})=V_{Max}(\{1\}) \\ \Leftrightarrow &  \text{  } \bigcap \mathfrak{B}=\{1\}.
	\end{align*}
\end{proof}

\begin{theorem}\label{1155}
	Let $\mathfrak{A}$ be a non-empty family of filters of a  residuated lattice $L$ such that $\mathfrak{A}\subseteq   Spec(L)$ and $\mathfrak{A}\cap E_L=\varnothing$, that is, $\mathfrak{A}$  is a non-empty family of prime filters that are not essential in $L$. Then the following are equivalent:
	\begin{enumerate}
		\item $\mathfrak{A}$ is a dense subset of $  Spec(L)$;
		\item   A filter $F$ of $L$ is essential in $L$ if and only if $F\not\subseteq P$ for each $P\in \mathfrak{A}$.
	\end{enumerate}
\end{theorem} 
\begin{proof} 
	$1\Rightarrow 2).$ Let  $F$ be an essential filter in $L$. Since every filter in $\mathfrak{A}$ is not essential in $L$,  we have $F\not\subseteq P$ for each $P\in \mathfrak{A}$ by  Remark \ref{eiclu}. Conversely, assume that $F\not\subseteq P$ for each $P\in \mathfrak{A}$. Choose $x\in F^*$. Then by  the definition, $[x)\cap F=\{1\}$. Hence, $[x)\cap F\subseteq P$ for each $P\in \mathfrak{A}$. Now since $F\not\subseteq P$ for each $P\in \mathfrak{A}$, we have $[x)\subseteq P$ for each $P\in \mathfrak{A}$, that is, $[x)\subseteq \bigcap_{P\in \mathfrak{A}}P$. By assumption and Corollary \ref{ddden}, $\bigcap_{P\in \mathfrak{A}}P=\{1\}$. Hence $x=1$, and so $F^*=\{1\}$. Thus $F$  is an essential filter in  $L$ by Proposition \ref{esse}.
	
	$2\Rightarrow 1).$ For each $P\in \mathfrak{A}$, $P^*\not=\{1\}$ by Proposition \ref{esse}. Hence, we can choose $1\not=x_P\in P^*$ for each $P\in \mathfrak{A}$. For each $P\in \mathfrak{A}$, if $x_P\in P$, then $x_P=x_P\vee x_P=1$, a contradiction. Hence, $x_P\not\in P$ for each $P\in \mathfrak{A}$. Consider the  $F:=[\{x_P\mid P\in \mathfrak{A}\})$.  For each $P\in \mathfrak{A}$ since $x_P\not\in P$, we have $F\not\subseteq P$ for each $P\in \mathfrak{A}$. Thus by assumption, the filter $F$ of $L$ is essential, that is, $F^*=\{1\} $ by Proposition \ref{esse}. Now let $x\in \bigcap \mathfrak{A}$. Hence $x\vee x_P=1$ for all $P\in \mathfrak{A}$. Thus $[x)\cap F=\{1\}$, and so we must have $[x)\subseteq F^*=\{1\}$. This shows that $x=1$, and hence $\bigcap \mathfrak{A}=\{1\}$. Therefore, $\mathfrak{A}$ is a dense subset of $  Spec(L)$ by Corollary \ref{ddden}.
\end{proof}   
\begin{proposition}\label{1166}
	Let  $P$  be a prime filter of a residuated lattice $L$. If  $P$ is not essential in $L$, then it is a minimal prime filter of $L$.
\end{proposition}
\begin{proof}
	Assume that  $P$ is a prime filter of a residuated lattice $L$ such that it  is not essential. By Proposition \ref{esse}, $P^*\not=\{1\}$. Choose $1\not=x\in P^*$. If  $x\in P$, then $x=x\vee x=1$, a contradiction. Hence,  $x\not\in P$. Now assume that $Q\subseteq P$ for some $Q\in  Spec(L)$. If $y\in P$, then $x\vee y=1$. Hence $x\vee y\in Q$. Since $x\not\in P$, we have $x\not\in Q$, and so $y\in P$. Thus $P\subseteq Q$. Therefore, $P=Q$ and so $P$ is a minimal prime filter of $L$.
\end{proof}

\begin{corollary}\label{pes}
	Let $P$ be a prime and finitely generated filter of a residuated lattice $L$. Then $P$ is not essential in $L$ if and only of it is a minimal prime filter of $L$.
\end{corollary}

\begin{proof}
	$\Rightarrow ).$ It follows from Proposition \ref{1166}.
	
	$\Leftarrow).$ Let $P$ be a minimal prime and  finitely generated filter of a residuated lattice $L$. Hence, there exist $x_1,..., x_n\in P$ such that $P=[x_1,..., x_n)$. By assumption and Theorem \ref{minp} for each $i\in\{1,..., n\}$, there exists $y_i\in L\setminus P$ such that $x_i\vee y_i=1$. Set $y:=y_1\vee\dots\vee y_n$. Since $P$ is a prime filter, we have $y\in L\setminus P$ and $x_i\vee y=1$ for each $i\in\{1,..., n\}$. Now if $a\in P$, then by Proposition \ref{2.555} there exist $m_1,..., m_n\in \mathbb{N}$ such that $x_1^{m_1}\odot\dots\odot x_n^{m_n}\leq a$. By Proposition \ref{2.222}, we have $x_i^{m_i}\vee y=1$ for each $i\in\{1,..., n\}$. Hence we have:\[1\odot\dots\odot 1= (x_1^{m_1}\vee y)\odot\dots\odot (x_n^{m_n}\vee y)) \leq (x_1^{m_1}\odot\dots\odot x_n^{m_n})\vee y \leq a\vee y.  \]
	Therefore,  $a\vee y=1$ and so $y\in P^*$. Now since $y\in L\setminus P$, we have $y\not=1$, and so $P^*\not=\{1\}$. Thus, $P$  is not essential in $L$ by Proposition \ref{esse}.
\end{proof}
\begin{theorem}\label{nm}
	Let $F$ be a filter of a  residuated lattice $L$ such that all minimal filter of $L$ is finitely generated (e.g., finite residuated latices). Then $F$ is an essential  filter of $L$ if and only if  $F\not\subseteq P$ for all minimal filter $P$ of $L$.
\end{theorem}

\begin{proof}
By assumption, every minimal prime  filter of $L$ is finitely generated. Hence $P\not\in E_L$ for all minimal prime  filter $P$ of $L$ by Corollary \ref{pes}. On the other hand, by Proposition \ref{2.700}, $\bigcap_{P \in  MinP(L)}P=\{1\}$. Hence the result follows from Corollary \ref{ddden} and Theorem \ref{1155}.
\end{proof}

\section{Independence of a family of filters}\label{2c}
We begin this section with the definition of an independent family of filters.
\begin{definition}
	A non-empty family $\{F_i\}_{i\in \Lambda}$ of filters of a residuated lattice  $L$ is called
	\textit{independent}, whenever  $F_j \cap(\bigvee_{j\not=i\in \Lambda}F_i)=\{1\}$ for
	each $j\in \Lambda$.
\end{definition}

\begin{definition}
	Let $F\in Filt(L)$ and $\{F_i\}_{i\in \Lambda}$ be a non-empty family of filters of a residuated lattice  $L$. Then we say $F$ \textit{is the direct sum of the family} $\{F_i\}_{i\in \Lambda}$ and denoted by $F=\bigoplus_{i\in \Lambda}F_i$, whenever $\{F_i\}_{i\in \Lambda}$ is an independent family of filters of $L$ and $F=\bigvee_{i\in \Lambda}F_i$. In this case, $F_i$ is called \textit{a direct summand of} $F$ for each $i\in \Lambda$.
\end{definition}

\begin{remark}
	If $\mathfrak{F}$ is an empty family of filters, then $\bigvee \mathfrak{F}=\{1\}$. Thus any family of filters with only one element is an independent family of filters. So each filter is a direct summand of itself.
\end{remark}
\begin{example}\label{exx}
Let $e\in  B(L)$ for a residuated lattice $L$. Then by Proposition \ref{2.555}, since $[e)\cap [e^*)=[e\vee e^*)=\{1\}$ and $[e)\cup [e^*)=[e\odot e^*)=[0)=L$, we have $L=[e)\oplus [e^*)$.
\end{example}
The following proposition gives a useful characterization for a two elements family of filters to be independent.
\begin{proposition}\label{indep1}
	Let $F,H\in  Filt(L)$ for a residuated lattice $L$. Then the following are equivalent:
	\begin{enumerate}
		\item $\{F, G\}$ is an independent family of filters of $L$;
		\item For $f, f'\in F$ and $g, g'\in G$, if $f\odot  g\leq f'\odot g'$, then $f\leq f'$  and $g\leq g'$;
	\item For $f, f'\in F$ and $g, g'\in G$, if $f\odot  g= f'\odot g'$, then $f= f'$  and $g= g'$.
	\end{enumerate}  
\end{proposition}
\begin{proof}
	$1\Rightarrow 2).$ Assume that  $f\odot  g\leq f'\odot g'$  for some $f, f'\in F$ and $g, g'\in G$. Then from $f\odot  g\leq f'\odot g'$, we have $f\odot  g\leq f'$ and $f\odot  g\leq  g'$. Thus $  g\leq f\to f'$ and $f\leq g\to g'$. Hence, $f\to f', g\to g'\in F\cap G=\{1\}$. Thus $f\to f'=1$ and $ g\to g'=1$, and so  $f\leq f'$  and $g\leq g'$.
	
	$2\Rightarrow 3).$ It is evident.
	
	$3\Rightarrow 1).$ Assume $x\in F\cap G$. Then since $x\odot 1= 1\odot x$ and $x, 1\in F\cap G$, we have $x=1$, and so $\{F, G\}$ is an independent family of filters of $L$. 
\end{proof}

\begin{corollary}\label{indep2}
	Let $F, F', H, H'\in  Filt(L)$ for a residuated lattice $L$ such that $F\subseteq F'$ and $H\subseteq H'$. Then if $F\oplus H=F'\oplus H'$, then $F= F'$ and $H= H'$.
\end{corollary}
\begin{proof}
	Let $a\in F'$ and $b\in H'$. Then $a\odot b\in F'\oplus H'$, and hence $a\odot b\in F\oplus H$. Thus by Proposition \ref{2.555}(5), there exist $f\in F$ and $h\in H$ such that $f\odot h\leq a\odot b$. Now since $a, f\in F'$ and $b, h\in H'$, we have $f\leq a$ and $h\leq b$ by Proposition \ref{indep1}.  From $f\in F$ and $h\in H$, we have  $a\in F$ and $b\in H$. Therefore, $F'\subseteq F$ and $H'\subseteq H$, and so  $F= F'$ and $H= H'$. 
\end{proof}

\begin{proposition}\label{11}
	Let $F,G\in  Filt(L)$ for a residuated lattice $L$. Then $L=F\oplus G$ if and only if there exists $e\in  B (L)$ such that $F=[e)$ and $G=[e^*)$. In particular, if $e, f\in B(L)$ and $L=[e)\oplus [f)$, then $f=e^*$. 
\end{proposition}
\begin{proof}
	$\Rightarrow).$  Since $L=F\oplus G$, we have $L=F\vee G$ and $H\cap G=\{1\}$. Thus, there exist $e\in F$ and $g\in G$ such that $e\odot g=0$. From $e\vee g\in F\cap G=\{1\}$, we have $e\vee g=1$. Hence by Proposition \ref{2.222}(5), $e\wedge g=e\odot g=0$. Hence $e\in B (L)$ and $e^*=g$.
	
	Now since $e\in F$ and $e^*\in G$, we have $[e)\subseteq F$ and $[e^*)\subseteq G$. If $x\in F$ and $y\in G$, then $0=e\odot e^*\leq x\odot y$. Hence by Proposition \ref{indep1}, we have $e\leq x$ and $e^*\leq y$. This shows that $x\in[e)$ and $y\in [e^*)$. Therefore, $F=[e)$ and $G=[e^*)$.

	$\Leftarrow).$ It is clear since $e\odot e^*=0$ and $e\vee e^*=1$. The particular case follows immediately from above argument, and the proof is complete.
\end{proof}
\begin{proposition}\label{1222}
	Let $F,G\in  Filt(L)$ for a residuated lattice $L$ such that $L=H\oplus G$. Then $F$ is maximal if and only if $G$ is simple. 
\end{proposition}
\begin{proof}
	By Proposition \ref{11}, there exists $e\in  B (L)$ such that $F=[e)$ and $G=[e^*)$. Hence the result follows from Proposition \ref{smb}.
\end{proof}
\begin{proposition}\label{2.7}
	Let $\{F_i\}_{i\in \Lambda}$ be a non-empty family of filters of a residuated lattice $L$. Then the following are equivalent:
	\begin{enumerate}
		\item $\{F_i\}_{i\in \Lambda}$ is an independent  family of filters of $L$;
		\item Every finite subfamily  $\{F_i\}_{i\in \Lambda'}$ of $\{F_i\}_{i\in \Lambda}$ is an independent  family of filters of $L$;
		
		\item If $A$ and $B$ are two non-empty subsets of $\Lambda$ such that $A\cap B=\varnothing$, then $(\bigvee_{i\in A}F_i)\cap(\bigvee_{j\in B}F_j)=\{1\}$;
	\item Let   $\{F_i\}_{i\in \Lambda'}$ be a  finite subfamily of $\{F_i\}_{i\in \Lambda}$ and $f_i, f'_i\in F_i$ for each $i\in \Lambda'$.  Then if $\odot_{i\in \Lambda'}f_i\leq \odot_{i\in \Lambda'}f'_i$, then $f_i\leq f'_i$ for each $i\in \Lambda'$;
	\item Let   $\{F_i\}_{i\in \Lambda'}$ be a  finite subfamily of $\{F_i\}_{i\in \Lambda}$ and $f_i, f'_i\in F_i$ for each $i\in \Lambda'$.  Then if $\odot_{i\in \Lambda'}f_i= \odot_{i\in \Lambda'}f'_i$, then $f_i= f'_i$ for each $i\in \Lambda'$.
	\end{enumerate}
\end{proposition}
\begin{proof}
	$1\Rightarrow 2).$ It is clear from the definition.
	
	$2\Rightarrow 3).$ Let $A$ and $B$ be two non-empty subsets of $\Lambda$ such that $A\cap B=\varnothing$. Choose $x\in(\bigvee_{i\in A}F_i)\cap(\bigvee_{j\in B}F_j)$. Then by Proposition \ref{2.555}(3), there exist finite sets $A'\subseteq A$ and $B'\subseteq B$ such that $x\in(\bigvee_{i\in A'}F_i)\cap(\bigvee_{j\in B'}F_j)$. Hence $x\in\bigvee_{i\in A'}F_i$ and $x\in\bigvee_{i\in B'}F_i$. Thus there exist $a_i\in F_i$ for each $i\in A'$ and  $b_j\in F_j$ for each $j\in B'$ such that $\odot_{i\in A'}a_i\leq x$ and $\odot_{j\in B'}b_j\leq x$. 
	For each $i\in A'$, $a_i\vee (\odot_{j\in B'}b_j)\in F_i\cap(\bigvee_{j\in B'}F_j) $. By assumption each $i\in A'$, we have $F_i\cap(\bigvee_{j\in B'}F_j)=\{1\} $. Thus for each $i\in A'$, $a_i\vee (\odot_{j\in B'}b_j)=1$. Hence by Proposition \ref{2.222}, we have:
	\begin{align*}
	1=\odot_{i\in A'}1= & \odot_{i\in A'}(a_i\vee (\odot_{j\in B'}b_j))\\ &\leq (\odot_{i\in A'}a_i)\vee (\odot_{j\in B'}b_j)\\ &\leq x.
	\end{align*}
	Hence $x=1$. Therefore,  $(\bigvee_{i\in A}F_i)\cap(\bigvee_{j\in B}F_j)=\{1\}$.
	
	$3\Rightarrow 1).$ It is clear from the definition, when for each $i\in \Lambda$, we set $A:=\{i\}$ and $B:=\{j\in \Lambda\mid i\not=j\}$.

$2\Leftrightarrow 4 \Leftrightarrow 5) $. It follows by using an inductive argument and Proposition \ref{indep1}.
\end{proof}

\begin{proposition}\label{2.10}
	Let $T$ be a simple filter of a residuated lattice $L$. Then we have the following:
	\begin{enumerate}
		\item There exists a prime filter $P$ of $L$ such that $\{T, P\}$ is a family of independent filters of $L$;
		\item If $Rad(L)=\{1\}$, then $T$ is a direct summand of $L$.
	\end{enumerate} 
\end{proposition}

\begin{proof}
	$1).$ By Proposition \ref{2.700}, $\bigcap_{P\in   Spec(L)}P=\{1\}$. Since $T\not=\{1\}$, we have $T\not\subseteq \bigcap_{P\in   Spec(L)}P$. Hence, there exists $P\in  Spec(L)$ such that $T\not\subseteq P$. Thus $T\cap P\not=T$. From the fact that $T$ is simple,  we have $T\cap P=\{1\}$. Therefore,  $\{T, P\}$ is a family of independent filters of $L$.	
	
	$2).$	Since $Rad(L)=\{1\}$ and $T\not=\{1\}$, we have $T\not\subseteq Rad(L)$. Hence, there exists $M\in Max(L)$ such that $T\not\subseteq M$. Thus $T\cap M\not=T$. From the fact that $T$ is simple, we have $T\cap M\subseteq T$ and $T\not\subseteq M$. Hence,  $T\cap M=\{1\}$ and $T\vee M=L$, that is, $L=T\oplus M$. Therefore, $T$ is a direct summand of $L$.
\end{proof}

We end this  section with two lemmas  that will be used in the following sections. 
\begin{lemma}\label{modul}
	Let $F, G, H\in  Filt(L)$ for a residuated lattice $L$ such that $F\subseteq G$. Then \[G\cap(F\vee H)=F\vee (G\cap H).\]
\end{lemma}
\begin{proof}
	Let $x\in G\cap(F\vee H)$. Then $x\in G$ and $x\in F\vee H$. Hence, there exist $f\in F$ and $h\in H$ such that $f\odot h\leq x$. Thus $h\leq f\to x$. Now since $h\leq f\to x$ and $x\leq f\to x$, we have $f\to x\in G\cap H$. Hence from $f\odot(f\to x)\leq x$, we have $x\in F\vee (G\cap H)$ by Proposition \ref{2.555}(5). Therefore, $G\cap(F\vee H)\subseteq F\vee (G\cap H)$. Conversely,  from $F\subseteq G$, we have $F\vee (G\cap H)\subseteq G$. Also, clearly $F\vee (G\cap H)\subseteq F\vee H$. Hence, $F\vee (G\cap H)\subseteq G\cap(F\vee H)$. Therefore $G\cap(F\vee H)=F\vee (G\cap H).$
\end{proof}
\begin{lemma}\label{sisc}
	Let $F\in Filt(L)$  for a residuated lattice $L$ such that for each $H\in Filt(L)$ with $H\subseteq F$,  there exists $G\in Filt(L)$  such that $F=H\oplus G$. Then for each $\{1\}\not=H\in Filt(L)$ with $H\subseteq F$, there exists $T\in   SimpF(L)$ such that $T\subseteq H$.
\end{lemma}
\begin{proof}
	Let $\{1\}\not=H\in Filt(L)$ with $H\subseteq F$. Choose $1\not=x\in H$ and set $\Sigma:=\{A\in Filt(L)\mid A\subsetneq [x)\}$. Since $\{1\}\in\Sigma$, we have $\Sigma\not=\varnothing$. Clearly $\Sigma$ is a partially ordered set by set inclusion and any chain of $\Sigma$ has an upper bound in $\Sigma$ (clearly, the union of any chain in $\Sigma$ is an upper bound for it). Thus, by
	Zorn's lemma, $\Sigma$ has a maximal element, say $B$. Since $B\subseteq F$, our assumption  follows that $F=B\oplus B'$ for some $B'\in Filt(L)$. So $[x)=[x)\cap F=[x)\cap(B\oplus B')=B\vee([x)\cap B')$ by Lemma \ref{modul}. Now if we show that $[x)\cap B'\in  SimpF(L)$, we are done since $[x)\cap B'\subseteq [x)\subseteq H$. 
	
	If $[x)\cap B'=\{1\}$, then $B=[x)$, a contradiction. Hence $[x)\cap B'\not=\{1\}$. Now let $G\subseteq [x)\cap B'$ for some $G\in Filt(L)$. We have two cases:
	
	$\mathbf{Case 1}$ If $G\subseteq B$, then $G\subseteq B\cap B'=\{1\}$, and so $G=\{1\}$ in this case.
	
	$\mathbf{Case 2}$ If $G\not\subseteq B$, then $B\subsetneq B\vee G\subseteq [x)$. By the selection of $B$, we must have $B\vee G= [x)$.  Choose $a\in [x)\cap B'$. Since $[x)\cap B'\subseteq [x)=B\vee G$,  there exist $b\in B$ and $g\in G$ such that $b\odot g\leq a$ by Proposition \ref{2.555}. Hence $b\leq g\to a$. From $a\leq g\to a$ and $b\leq g\to a$, we have $g\to a\in B\cap  [x)\cap B'=\{1\}$. Thus $g\to a=1$, and so $g\leq a$. Now since $g\in G$, we have $a\in G$. Therefore, $[x)\cap B'\subseteq G$, and so $[x)\cap B'= G$ in this case. 
	
	Therefore,  $[x)\cap B'\in  SimpF(L)$.
\end{proof}

\section{The Socle of a filter}\label{3c}
In this section, the socle of a filter is introduced, and then we consider its relation with essential filters.
\begin{definition}
	The \textit{socle} of  a filter $F$ of a residuated lattice $L$, denoted by $ Soc((F)$, is equal to $\bigvee_{ T\in  SimpF(L), T\subseteq F}T$ 
	if the set $\{T\mid T\in  SimpF(L), T\subseteq F\}$ is non-empty, and is equal to $\{1\}$ if  $\{T\mid  T\in  SimpF(L), T\subseteq F\}=\varnothing$.
\end{definition}

The next theorem states the relation between the socle of a filter and essential filters.

\begin{theorem}\label{3.3}
	Let $F\not=\{1\}$ be a filter of a residuated lattice $L$. Then \[ Soc((F)=\bigcap E_F.\]
\end{theorem}
\begin{proof}
	If $\{T\in  SimpF(L)\mid T\subseteq F\}=\varnothing$, then $ Soc((F)=\{1\}$ and we have $ Soc((F)\subseteq \bigcap E_F$ in this case. 
	
	Now assume that  $\{T\in  SimpF(L)\mid T\subseteq F\}\not=\varnothing$. If $A\in \{T\in  SimpF(L)\mid T\subseteq F\}$ and $B\in E_F$, by the definition of $E_F$, if $A\cap B=\{1\}$, then $A\cap F=\{1\}$. Now since $A\subseteq F$, we have $A=\{1\}$, a contradiction. Hence for each $A\in \{T\in  SimpF(L)\mid T\subseteq F\}$ and $B\in E_F$, we have  $A\cap B\not=\{1\}$, and since  each $A\in \{T\in  SimpF(L)\mid T\subseteq F\}$ is a simple filter of $L$,  for each $A\in \{T\in  SimpF(L)\mid T\subseteq F\}$ and $B\in E_F$, we have  $A\cap B=A$, and so $A\subseteq B$. Thus $ Soc((F)\subseteq\bigcap E_F$.
	
	Let $H\in Filt(L)$ such that $H\subseteq \bigcap E_F$. Now since $F\in E_F$, we conclude that $H\subseteq F$, hence that $H\vee(H^*\cap F)\in E_F$ by Lemma \ref{efe}, and finally that $H\subseteq \bigcap E_F\subseteq H\vee(H^*\cap F)\in E_F$. Thus by Lemma \ref{modul}, we have \[\bigcap E_F=\bigcap E_F\cap(H\vee(H^*\cap F))=H\vee((\bigcap E_F)\cap H^*\cap F).\] 

From $H\cap ((\bigcap E_F)\cap H^*\cap F)\subseteq H\cap H^*=\{1\}$, we have $H\cap ((\bigcap E_F)\cap H^*\cap F)=\{1\}$. Hence \[\bigcap E_F=H\oplus((\bigcap E_F)\cap H^*\cap F).\] Thus the condition of  Lemma \ref{sisc}	for the filter $\bigcap E_F$ is satisfied. From $ Soc((F)\subseteq\bigcap E_F$, we have  $\bigcap E_F= Soc((F)\oplus((\bigcap E_F)\cap  Soc((F)^*\cap F)$.

	Set $A:=(\bigcap E_F)\cap  Soc((F)^*\cap F$. If $A=\{1\}$, then $ Soc((F)=\bigcap E_F$ and we are done. Now assume that $A\not=\{1\}$. Then by Proposition \ref{sisc}, there exists $T\in  SimpF(L)$ such that $T\subseteq A$. Hence, \[\{1\}\not=T\subseteq  Soc((F)\cap((\bigcap E_F)\cap  Soc((F)^*\cap F)=\{1\},\] 
	
	which is impossible. 
\end{proof}
\begin{theorem}\label{3.4}
	Let $F$ be a filter of a residuated lattice $L$. Then we have the following:
	\begin{enumerate}
		\item $ Soc((F)=F\cap Soc((L)$;	
		\item $ Soc(( Soc((F))= Soc((F)$.
	\end{enumerate}
\end{theorem}
\begin{proof}
	$1)$ Since $ Soc((F)\subseteq F$  and $ Soc((F)\subseteq  Soc((L)$, we have $ Soc((F)\subseteq F\cap Soc((L)$. Now let $x\in F\cap Soc((L)$. Hence $x\in  Soc((L)$. Then by Proposition \ref{2.555}(3), there exist $T_1,..., T_n\in  SimpF(L)$ and $x_1\in T_1,...,x_n\in T_n$ such that $x_1\odot\dots\odot x_n\leq x$. Set $I:=\{i\mid i\in\{1,...,n\}\text{  and  } x\vee x_i=1\}.$ We consider the following cases:
	
	$\mathbf{Case\text{ } 1.}$ If $I=\varnothing$, then $x\vee x_i\not=1$ for all $i\in \{1,..., n\}$. Hence, $1\not=x\vee x_i\in F\cap T_i$, and so $F\cap T_i\not=\{1\}$. From $F\cap T_i\subseteq T_i$ and the fact that $T_i$ is a simple filter, we have $F\cap T_i=T_i$. Thus $T_i\subseteq F$ for all $i\in \{1,..., n\}$, that is,  $T_i\subseteq  Soc((F)$ for all $i\in \{1,..., n\}$. Therefore, $x\in  Soc((F)$ in this case.
	
	$\mathbf{Case\text{ } 2.}$ If $I=\{1,..., n\}$, then $x\vee x_i=1$ for all $i\in \{1,..., n\}$. Then by Proposition \ref{2.555}(6), we have 
	\[	x=x\vee(x_1\odot\dots\odot x_n)\geqslant(x\vee x_1)\odot\dots\odot(x\vee x_n)=1.\]
	Thus $x=1$, and so $x\in  Soc((F)$ in this case.
	
	$\mathbf{Case\text{ } 3.}$ If $I\not=\varnothing$ and $I\not=\{1,..., n\}$. Set $J:=\{1,..., n\}\setminus
	I$. Hence by Proposition \ref{2.555}(6), we have:
	\begin{align*}
	(\odot_{j\in J}(x\vee x_j))= & (\odot_{i\in I}(x\vee x_i))\odot(\odot_{j\in J}(x\vee x_j))\\&=  (x\vee x_1)\odot\dots\odot(x\vee x_n) \\&\leq  x\vee(x_1\odot\dots\odot x_n)=x.\nonumber
	\end{align*} Hence
	\begin{align}\label{sss}
	(\odot_{j\in J}(x\vee x_j))\leq  x\vee(x_1\odot\dots\odot x_n)=x.\tag{*}
	\end{align}
	Now since $x\vee x_j\not=1$ for each $j\in J$. Hence, $1\not=x\vee x_j\in F\cap T_j$, and so $F\cap T_j\not=\{1\}$ for each $j\in J$. From $F\cap T_j\subseteq T_j$ and the fact that $T_j$ is a simple filter, we have $F\cap T_j=T_j$. Thus $T_j\subseteq F$ for each $j\in J$, that is, $x\vee x_j\in  Soc((F)$ for each $j\in J$. Thus $(\odot_{j\in J}(x\vee x_j))\in  Soc((F)$. Therefore, $x\in  Soc((F)$ by $(*)$ in this case.
	
	$2)$ By part $(1)$, Since $ Soc((F)\subseteq Soc((L)$ we have \[ Soc(( Soc((F))= Soc((F)\cap Soc((L)= Soc((F).\]
\end{proof}

\begin{theorem}\label{miso}
	Let $M$ be a maximal filter of a residuated lattice $L$ with $Rad(L)=\{1\}$.
	Then the following are equivalent:
	
	\begin{enumerate}
		\item $M$ is an isolated point of $Max(L)$;
		\item There exists $e\in B (L)$ such that $M=[e^*)$;
		\item $M$ is a direct summand of $L$.
	\end{enumerate}
\end{theorem}
\begin{proof}
	$1\Rightarrow 2)$ Let $M$ be an isolated point of $Max(L)$. Hence there exists $a\in L$ such that $\{M\}=U_{Max}(a)$, that is, $a\in(\bigcap_{M\not= N\in Max(L)}N)\setminus M$. Since $a\not\in M$, there exists $n\in\mathbb{N}$ such that $(a^n)^*\in M$ by Theorem \ref{tpd}. Set $e:=a^n$. Hence, $e\in(\bigcap_{M\not= N\in Max(L)}N)\setminus M$ and $e^*\in M$. Thus $e\vee e^*\in \bigcap_{ N\in Max(L)}N=Rad(L)=\{1\}$, that is, $e\vee e^*=1$, and so $e\in  B(L)$ by Proposition \ref{2.222}(13). 
	Now since $e^*\in M$, we have $[e^*)\subseteq M$. If $x\in M$, then $x\vee e\in Rad(L)=\{1\}$. Hence $x\vee e=1$. Thus by Proposition \ref{2.222}, we have $e^*\to x=e^{**}\vee x=e\vee x=1$. Hence, $e^*\leq x$, and so $x\in[e^*)$. Therefore, $M=[e^*)$.
	
	$2\Rightarrow 1)$ It is clearly that if there exists $e\in B(L)$ such that $M=[e^*)$, then we have $\{M\}=U_{Max}(e)$, and so $M$ is an isolated point of $Max(L)$.
	
	$2\Rightarrow 3)$ If there exists $e\in B (L)$ such that $M=[e^*)$, then by Example \ref{exx}, we have $L=[e)\oplus[e^*)$, and so $M$ is a direct summand of $L$.
	
	$3\Rightarrow 2)$ It follows from Proposition \ref{11}.
\end{proof}

\begin{theorem}\label{3.6}
	Let $L$ be a residuated lattice with $ Rad(L)=\{1\}$. Then 
	\[ Soc((L)=\{a\in L\mid U_{max}(a) \text{  }\text{is a finite set}\}.\]
\end{theorem}
\begin{proof}
	If $ Soc((L)=\{1\}$, then $ Soc((L)\subseteq\{a\in L\mid U_{max}(a) \text{  }\text{is a finite set}\}$. Now, assume that $ Soc((L)\not=\{1\}$ and 
 $1\not=a\in Soc((L)=\bigvee_{T\in   SimpF(L)}T$, the there exist $T_1,..., T_n\in   SimpF(L)$ and $a_1\in T_1$,..., $a_n\in T_n$ such that $a_1\odot\dots\odot a_n\leq a$. Hence, $a\in T_1\vee\dots\vee T_n$ and so by Proposition \ref{smss}, we have $U_{max}(a)\subseteq U_{max}(T_1\vee\dots\vee T_n)$. By Proposition \ref{smss}, $U_{max}(T_1\vee\dots\vee T_n)=U_{max}(T_1)\cup\dots\cup U_{max}( T_n)$. From  Proposition \ref{133}, $U_{max}(T_i)$ is a finite set  for each $i$, hence $U_{max}(a)$ is a finite set, and finally that $ Soc((L)\subseteq\{a\in L\mid U_{max}(a) \text{  }\text{is a finite set}\}$. 
	
	Conversely, assume that $U_{max}(a)$ is a finite set for $a\in L$, say $U_{max}(a)=\{M_1,...M_n\}$. Hence for each $N\in Max(L)\setminus U_{max}(a)$, we have $a\in N$. If $n=1$, then $M_1$ is an isolated point of $Max(L)$, and so there exists $e\in B (L)$ such that $M=[e)$ by Theorem \ref{miso}. Now since $a\vee e\in  Rad(L)=\{1\}$, we have $a\vee e=1$. Thus  by Proposition \ref{2.222},  $e^*\leq a$ and hence $a\in [e^*)$. For Proposition \ref{smb} and the fact that $[e)$ is a maximal filter, $[e^*)$ is a simple filter and so $a\in Soc((L)$ in this case. Now assume that $2\leq n$.  Let $j\in\{1,..., n\}$ be fixed. Since $M_1,..., M_n$ are maximal filters, there exists $a_i\in M_i\setminus M_j $ for each $i\not= j$ and $i\in \{1,..., n\}$. Set $b_j:=a\vee a_1\vee\dots \vee \hat{a_{j}}\vee\dots\vee a_n$, where the hat means that term is omitted. Clearly $b_j\in (\bigcap_{M_j\not=N\in Max(L)}N)\setminus M_j$. Hence $U_{max}(b_j)=\{M_j\}$ for each $j\in\{1,..., n\}$. Thus for each $j\in\{1,..., n\}$,  $M_j$ is an isolated point of $Max(L)$, and so there exists $e_j\in B (L)$ such that $M_j=[e_j)$ by Theorem \ref{miso}. Now since $a\vee e_1\vee\dots\vee e_n\in  Rad(L)=\{1\}$, we have $a\vee e_1\vee\dots\vee e_n=1$. Thus by Proposition \ref{2.222}, $ (e_1\vee\dots\vee e_n)^*\leq a$, that is, $ (e_1^*\odot\dots\odot e_n^*)\leq a$ and hence $a\in [e_1^*)\vee\dots\vee[e_n^*)$. From Proposition \ref{smb} and the fact that $[e_j)$ is a maximal filter, $[e_j^*)$ is a simple filter and so $a\in Soc((L)$ in this case. 
\end{proof}

\begin{theorem}\label{3.7}  
	Let $L$ be a residuated lattice. Then the following are equivalent:
	\begin{enumerate}
		\item $ Soc((L)\in E_L$, that is, $ Soc((L) $ is an essential filter of $L$;
		\item For each $\{1\}\not=F\in  Filt(L)$, we have $ Soc((F)\not=\{1\}$;
		\item  Every non-trivial filter of $L$ contains a simple filter;
		\item The intersection of any non-empty subfamily of essential filters is an essential filter, that is, the intersection of any non-empty subfamily of $E_L$ is a member of $E_L$.
	\end{enumerate}

\end{theorem}

\begin{proof}
	$1\Rightarrow 2).$ By Theorem \ref{3.4},  $ Soc((F)=F\cap Soc((L)$. If $ Soc((F)=\{1\}$, then $F\cap Soc((L)=\{1\}$. Now by assumption since $ Soc((L) $ is an essential filter of $L$,  we have $F=\{1\}$.
	
	$2\Rightarrow 3).$ It is clear by the definition of the socle of a filter.
	
	$3\Rightarrow 1).$ If $F$ is a non-trivial filter of $L$, then $F$ contains a simple filter by our assumption. Hence, $F\cap Soc((L)\not=\{1\}$, that is, $( Soc((L))^*=\{1\}$, and so $ Soc((L)\in E_L$ by Proposition \ref{esse}.
	
	$1\Rightarrow 4).$ Let $ Soc((L)\in E_L$. If $\mathfrak{A}$ is a non-empty subfamily of $E_L$, then $ Soc((L)\subseteq\bigcap \mathfrak{A}$. Thus, $\bigcap \mathfrak{A}$  is an essential filter of $L$ by Remark \ref{eiclu}.
	
	$4\Rightarrow 1).$ By Theorem \ref{3.3}, $ Soc((L)=\bigcap E_L$. Hence $ Soc((L)\in E_L$ by our assumption.
\end{proof} 
Recall that a residuated $L$ is called \textit{Artinian}, whenever every decreasing chain of filters of $L$  is stationary, see [1].
\begin{corollary}\label{3.8}
	$ Soc((L)\in E_L$ 	for any non-trivial Artinian (e.g., finite) residuated lattice $L$.
\end{corollary}
\begin{proof}
	It is clear that any non-trivial filter of a non-trivial Artinian residuated lattice contains a simple filter. Hence the result follows from Theorem \ref{3.7}.
\end{proof}

For a residuated lattice $L$, we set $M_0(L):=\{M\in Max(L)\mid M \text{ }\text{is an isolated point of}\text{ } Max(L)\}$.

\begin{theorem}\label{3.9}
	Let $L$ be a residuated lattice with $ Rad(L)=\{1\}$. Then  $M_0(L)$ is a dense subset of $Max(L)$ is equivalent to the (equivalent) conditions in Theorem \ref{3.7}. 
	
\end{theorem} 

\begin{proof}
 Assume that $M_0(L)$ is a dense subset of $Max(L)$, hence that  $\bigcap M_0(L)=\{1\}$ by Corollary \ref{dddenm}. By Corollary \ref{sr1} for each $T\in  SimpF(L)$, there exists $e_T\in B (L)$ such that $T=[e_T)$. By Proposition  \ref{smb}, $\{[e^*_T)\mid T\in  SimpF(L)\}\subseteq Max(L)$. Now let $x\in( Soc((L))^*$. Then $x\vee e_T=1$ for each $T\in  SimpF(L)$. Thus $e_T^*\leq x$ for each $T\in  SimpF(L)$, that is, $x\in[e_T^*)$ for each $T\in  SimpF(L)$. Thus $x\in \bigcap M_0(L)=\{1\}.$ This shows that $x=1$, and so $ Soc((L)\in E_L$ by Proposition \ref{esse}.

	Conversely, assume that  $ Soc((L)\in E_L$ and $x\in \bigcap M_0(L)$. Now if $y\in  Soc((L)$, then by Proposition \ref{2.555} and Corollary  \ref{sr1},  there there exist $e_1,..., e_n\in B (L)$ such that $[e_1),..., [e_n)\in  SimpF(L)$ and $e_1\odot\dots\odot e_n\leq y$. 
	By  Proposition \ref{smb} and Theorem \ref{miso} for each $i$, $[e_i^*)$ is a maximal filter and is an isolated point of $Max(L)$. Hence $x\in [e_i^*)$, that is, $e_i^*\leq x$. So $x\vee e_i=1$ for each $i\in\{1,..., n\}$. By Proposition \ref{2.222}, we have \[1=(x\vee e_1)\odot\dots\odot(x\vee x_n)\leq x\vee(e_1\odot\dots\odot e_n)\leq x\vee y, \]
	and hence $x\vee y=1$. Therefore,  $x\vee y=1$ for each $y\in  Soc((L)$. Hence by our assumption  and Proposition \ref{esse}, we have $x=1$, and so $\bigcap M_0(L)=\{1\}$. Thus $M_0(L)$ is a dense subset of $Max(L)$ by Corollary \ref{dddenm}.
	
\end{proof}

\begin{theorem}\label{3.10}
	Let $L$ be a residuated lattice with $ Rad(L)=\{1\}$. Then $ Soc((L)$ is a principal filter of $L$ if and only if $M_0(L)$ is a finite subset of $Max(L)$.
\end{theorem}
\begin{proof}
	$\Rightarrow ).$ Assume that $ Soc((L)$ a principal filter of $L$. Hence, there exists $a\in L$ such that $ Soc((L)=[a)$. By Theorem \ref{miso} and our assumption, every member of $M_0(L)$ is generated by a Boolean element, hence assume that $M_0=\{[e_i)\}_{i\in\Lambda}$, where $\{e_i\}_{i\in\Lambda}\subseteq  B (L)$. By definition, $ Soc((L)=\bigvee_{T\in  SimpF(L)}T$. From Corollary  \ref{sr1} since  $a\in Soc((L)$, there exist $T_1,..., T_n\in  SimpF(L)$ and $f_1,..., f_n\in B (L)$ such that $T_j=[f_j)$ for each $j\in\{1,..., n\}$ and $f_1\odot\dots\odot f_n\leq a$. Hence for each $i\in \Lambda$ since $[e_i)$ is a maximal filter, $[e_i^*)$ is a minimal filter, and so $e_i^*\in Soc((L)$ for each $i\in \Lambda$. Now assume that $i\in \Lambda$, then $e_i^*\in Soc((L)=[a)$. Thus there exists $m\in\mathbf{N}$ such that $a^m\leq e_i^*$. Hence,
	\begin{align*}
	f_1\wedge\dots\wedge f_n= & f_1\odot\dots\odot f_n=(f_1\odot\dots\odot f_n)^m\leq a^m\leq e_i^*\\& \Rightarrow f_1\wedge\dots\wedge f_n\leq e_i^*\\ & \Rightarrow e_i\leq (f_1\wedge\dots\wedge f_n)^*
	\\ & \Rightarrow e_i\leq f_1^*\vee\dots\vee f_n^*
	\\ & \Rightarrow [f_1^*\vee\dots\vee f_n^*)\subseteq [e_i)
	\\ & \Rightarrow [f_1^*)\cap\dots\cap [f_n^*)\subseteq [e_i).
	\end{align*}
	Now since 	$[e_i)$ is maximal, it is prime, and so there exists $j\in\{1,..., n\}$ such that $[f_j^*)\subseteq [e_i)$. Now since $[f_j)$ is simple, by Proposition \ref{smb},  $[f_j^*)$ is maximal. Also since $[e_i)$ is a maximal filter, we have $[f_j^*)= [e_i)$. So $e_i=f_j^*$ by Proposition \ref{2.555}. This shows that $M_0(L)$ is a finite subset of $Max(L)$.
	
	$\Leftarrow ).$ Assume that $M_0(L)$ is a finite subset of $Max(L)$. Since $ Rad(L)=\{1\}$, every simple filter of $L$ is of the form $[e)$ for some $e\in B (L)$ such that $[e^*) \in M_0(L)$ by Proposition \ref{smb} and Theorem \ref{miso}. Now since  $M_0(L)$ is a finite, $  SimpF(L)$ is also finite. Assume that $  SimpF(L)=\{[e_1),..., [e_n)\}$. Then \[ Soc((L)=[e_1)\vee\dots\vee [e_n)=[e_1\odot\dots\odot e_n).\] Hence $ Soc((L)$ a principal filter of $L$. 
\end{proof}

\section{Semi-simple filters}\label{4c}
In this section, we consider semi-simple filters of a residuated lattice. We begin with the following definition. 

\begin{definition}
	A filter $F$ of a residuated lattice $L$ is called \textit{semi-simple}, whenever $F\not=\{1\}$ and $F=\bigvee_{i\in \Lambda}T_i$ for some non-empty family  $\{T_i\}_{i\in \Lambda}$ of simple filters.
\end{definition}

In the next theorem, we state some equivalent conditions for a filter $F$ of a residuated lattice $L$  being semi-simple. 

\begin{theorem}\label{4.2}
	Let $F$ be a filter of a residuated lattice $L$. Then the following are equivalent:
	\begin{enumerate}
		\item $F=\bigvee_{i\in \Lambda}T_i$ for some non-empty family  $\{T_i\}_{i\in \Lambda}$ of simple filters, that is, $F$ is semi-simple;
		
		\item $F=\bigoplus_{j\in \Omega}T_j$ for some non-empty family  $\{T_i\}_{j\in \Omega}$ of simple filters;
		\item If $H\in Filt(L)$ and $H\subseteq F$, then there exists $G\in Filt(L)$  such that $F=H\oplus G$;
	\item If $H\in Filt(L)$ and $H\subseteq F$, then  $F=H\oplus (H^*\cap F)$;
		\item $E_F=\{F\}$.
	\end{enumerate} 
\end{theorem}
\begin{proof}
	$1\Rightarrow 2)$ Assume that $F=\bigvee_{i\in \Lambda}T_i$ for some non-empty family  $\{T_i\}_{i\in \Lambda}$ of simple filters. Set $\Sigma:=\{K\subseteq \Lambda\mid \bigvee_{k\in K}T_k=\bigoplus_{k\in K}T_k\}$. Clearly $\Sigma$ is a partially ordered set by set inclusion, and for each $i\in \Lambda$, $\{i\}\in\Sigma$, and so $\Sigma\not=\varnothing$. Let $\{\Gamma_\alpha\}$ be a chain in $\Sigma$. Set $\Gamma:=\bigcup_{\alpha}\Gamma_\alpha$. For $\beta\in \Gamma$ if $x\in T_\beta\bigcap(\bigvee_{\beta\not=\alpha\in \Gamma}T_\alpha)$, then there exists $\alpha_1,..., \alpha_n\in \Gamma$ such that $x\in T_\beta\bigcap(T_{\alpha_1}\vee\dots\vee T_{\alpha_n} )$. Hence there exists $\gamma$ such that $\beta, \alpha_1,..., \alpha_n\in \Gamma_\gamma$, and so $x\in T_\beta\bigcap(T_{\alpha_1}\vee\dots\vee T_{\alpha_n} )=\{1\}$. Thus $\Gamma\
\in \Sigma$. Hence every chain in $\Sigma$ has an upper bound in $\Sigma$. Thus, by Zorn's Lemma, $\Sigma$ has a maximal element, say $\Omega$. Set $G:=\bigvee_{a\in \Omega}T_a=\bigoplus_{a\in \Omega}T_a$. If $F\not=G$, then there exists $i\in \Lambda$ such that $T_i\not\subseteq G$. Since  $T_i$ is simple and $T_i\not\subseteq G$, we have $T_i\bigcap(\bigvee_{a\in \Omega}T_a)=\{1\}$. If $x\in T_{a_*}\bigcap(T_i\vee(\bigvee_{a_*\not=a\in \Omega}T_a))$ for some $a_*\in \Omega$, then $x\in T_i\vee(\bigvee_{a_*\not=a\in \Omega}T_a)$ and so there exist $t\in T_i$ and $l\in \bigvee_{a_*\not=a\in \Omega}T_a$ such that $t\odot l\leq x$, hence $t\leq l\to x$. So $l\to x\in T_i\cap (\bigvee_{a_*\not=a\in \Omega}T_a)\subseteq T_i\bigcap(\bigvee_{a\in \Omega}T_a)=\{1\} $. This shows that $l\to x=1$, and so $l\leq x$. Thus $ x\in T_{a_*}\cap(\bigvee_{a_*\not=a\in \Omega}T_a)=\{1\} $. Hence $x=1$ and so $\Omega\cup \{i\}\in \Sigma$, a contradiction. Therefore $G=F$, and so $F=\bigoplus_{j\in \Omega}T_j$, where  $\{T_j\}_{j\in \Omega}\subseteq   SimpF(L)$.
	
	$2\Rightarrow 1)$ It is clear.
	
	$2\Rightarrow 3)$ Let $F=\bigoplus_{j\in \Omega}T_j$ for some non-empty family  $\{T_i\}_{j\in \Omega}$ of simple filters. Choose  $H\in Filt(L)$ and $H\subseteq F$. If $H=F$, then we have nothing to prove. So assume that $H\not= F$. Set $\Sigma:=\{K\subseteq \Omega\mid H\cap (\bigvee_{k\in K}T_k)=\{1\}\}$. Clearly $\Gamma$ is a partially ordered set by set inclusion and since $H\not= F$, there exists $j\in \Omega$ such that $T_j\not\subseteq H$, and so $\{j\}\in\Gamma$, hence that $\Sigma\not=\varnothing$.

	Let $\{\Gamma_\alpha\}$ be a chain in $\Sigma$. Set $\Gamma:=\bigcup_{\alpha}\Gamma_\alpha$. If $x\in H\cap(\bigvee_{j\in \Gamma}T_j)$, then there exist $j_1,..., j_n\in \Gamma$ such that $x\in H\cap(T_{j_1}\vee\dots\vee T_{j_n} )$. Hence there exists $j_*\in \Gamma$ such that $ j_1,..., j_n\in \Gamma_{j_*}$, and so $x\in H\cap(\bigvee_{j\in \Gamma_{j_*}}T_j)=\{1\}$. Thus $x=1$, and hence  every chain in $\Sigma$ has an upper bound in $\Sigma$. Thus, by
	Zorn's lemma, $\Sigma$ has a maximal element, say $A$. Now set $G:=\bigvee_{a\in A}T_a$. Hence $H\cap G=\{1\}$. By a similar argument in $(1\Rightarrow 2)$ $F=H\vee G$. Therefore, $F=H\oplus G$.

	$3\Rightarrow 1)$ By Lemme \ref{sisc}, for each $\{1\}\not=H\in Filt(L)$ with $H\subseteq F$ there exists $T\in   SimpF(L)$ such that $T\subseteq H$. Set $G:=\bigvee_{T\in  SimpF(L), T\subseteq F}T$. By assumption, $F=G\oplus G'$ for some $G'\in  Filt(L)$. If $G'\not=\{1\}$, then by Lemma \ref{sisc} there exists $T\in   SimpF(L)$ such that $T\subseteq G$. Hence $T\subseteq G\cap G'=\{1\}$, and so $T=\{1\}$, a contradiction. Thus, $G'=\{1\}$, and so $F=\bigvee_{T\in  SimpF(L), T\subseteq F}T$, and we are done.

	$3\Rightarrow 5)$ Clearly $F\in E_F$. Let $H\in E_F$. By assumption then there exists $G\in Filt(L)$  such that $F=H\oplus G$. Since $H\cap G=\{1\}$, we $F\cap G=\{1\}$. From $G\subseteq F$, we   have $G=\{1\}$, and so $F=H$. Hence $E_F=\{F\}$.
	
	$5\Rightarrow 4)$ Let  $H\in Filt(L)$ and $H\subseteq F$. By Lemma \ref{efe}, $H\vee(H^*\cap F)\in E_F$. So $H\vee(H^*\cap F)=F$. Now since $H\cap(H^*\cap F)=\{1\}$. we have $F=H\oplus(H^*\cap F)$, and we are done.
	
$4\Rightarrow 3)$	It is clear.
\end{proof} 

\begin{corollary}\label{4.3}
	Let $F, G\in Filt(L)$ such that $G\subseteq H$. Then if $F$ is semi-simple, then $G$ is also semi-simple.
\end{corollary}
\begin{proof}
	Assume that $H\in Filt(L)$ such that $H\subseteq G$. Then  $H\subseteq F$. Hence by Theorem \ref{4.2}, there exists $H'\in  Filt(L)$ such that $F=H\oplus H'$. Thus we have 
	$H\subseteq G\subseteq H\oplus H'=H\vee H'$, and so by Lemma \ref{modul} $G=G\cap(H\vee H')=H\vee(G\cap H')$. Now since $H\cap(G\cap H')=\{1\}$, we have $G=H\oplus(G\cap H')$. Therefore, $G$ is  semi-simple by Theorem \ref{4.2}.
\end{proof}

In the following theorem, we state some equivalent conditions for the filter $L$ of a residuated lattice $L$  being semi-simple.

\begin{theorem}\label{4.4}
	Let $L$ be a non-trivial  residuated lattice $L$. Then the followings are equivalent:
	\begin{enumerate}
		\item $L=\bigvee_{i\in \Lambda}T_i$ for some non-empty family  $\{T_i\}_{i\in \Lambda}$ of simple filters, that is, $L$ is semi-simple as a filter of $L$;
		
		\item $L=\bigoplus_{j\in \Omega}T_j$ for some non-empty family  $\{T_i\}_{i\in \Omega}$ of simple filters;
		\item If $H\in Filt(L)$, then there exists $G\in Filt(L)$  such that $L=H\oplus G$;
	\item If $H\in Filt(L)$, then  $L=H\oplus H^*$;
	\item Every filter of $L$ can be generated by a Boolean element;
		\item $E_L=\{L\}$;
		
		\item $L=T_1\oplus\dots\oplus T_n$ for some $T_1,..., T_n\in   SimpF(L)$;
		
		\item  $L$ is hyperarchimedean and semi-local;
		\item  $ Rad(L)=\{1\}$ and  $L$ semi-local;
		\item Every proper filter of $L$ can be written uniquely, up to order, as a finite intersection of distinct maximal filters.
		\item The filter $\{1\}$ can be written  as a  finite intersection of distinct maximal filters;
		\item Every non-trivial filter of $L$ is semi-simple.
	\end{enumerate} 
\end{theorem}
\begin{proof}

	By Theorem \ref{4.2}, we have  $1\Leftrightarrow 2\Leftrightarrow 3\Leftrightarrow 4\Leftrightarrow 6$.

$5\Leftrightarrow 6)$ It follows frome Example \ref{exx} and Proposition \ref{11}.

	$2\Rightarrow 7)$ Assume that $L=\bigoplus_{j\in \Omega}T_j$ for some non-empty family  $\{T_i\}_{i\in \Omega}$ of simple filters. Since $0\in L$, there exist $j_1,..., j_n\in \Omega$ such that $0\in T_{j_1}\vee\dots\vee T_{j_n}$, that is, $L=T_{j_1}\vee\dots\vee T_{j_n}$. Now since $\{T_j\}_{j\in \Omega}$ is an independent family of filters of $L$, $\{T_{j_1},..., T_{j_n}\}$ is also an independent family of filters of $L$ by Proposition \ref{2.7}. Hence $L=T_{j_1}\oplus\dots\oplus T_{j_n}$. 
	
	$7\Rightarrow 2)$ It is clear. 
	
	$3\Rightarrow 8$ Let $P\in  Spec(L)$. Choose $M\in Max(L)$ such that $P\subseteq M$. By assumption, there exists $G\in  Filt(L)$ such that $L=M\oplus G$. Now if $P\not=M$, then there exists $x\in M\setminus P$. If $y\in G$, then $x\vee y\in M\cap G=\{1\}$, and so $x\vee y=1\in P$. Since $x\not\in P$, we have $y\in P$. Hence  $G\subseteq P$. From $G\cap M=\{1\}$ and $G\subseteq P\subseteq M$, we have $G=\{1\}$. Thus $L=M\vee G=M\vee \{1\}=M$, a contradiction. Hence, $P=M$, and so $L$ is hyperarchimedean by Theorem \ref{hyp}.
	
	Now since $(3\Leftrightarrow 7)$,  $L=T_1\oplus\dots\oplus T_n$ for some $T_1,..., T_n\in   SimpF(L)$. By Theorem \ref{hyp1} since $L$ is hyperarchimedean, we have $ Rad(L)=\{1\}$, and so by Corollary \ref{sr1} there exist $e_1,..., e_n\in  B (L)$ such that $T_i=[e_i)$ for each $i\in\{1,..., n\}$. Hence for each $i\in\{1,..., n\}$ we have 
	\begin{align*}
	L=T_1\oplus\dots\oplus T_n= & [e_1)\oplus\dots\oplus [e_n)\\=&[e_i)\oplus[e_1)\oplus\dots\oplus\hat{[e_i)}\oplus\dots\oplus[e_n)
	\\=&[e_i)\oplus[e_1)\vee\dots\vee\hat{[e_i)}\vee\dots\vee  [e_n)\\=& [e_i)\oplus[e_i\odot\dots\odot\hat{e_i}\odot\dots\odot  e_n),
	\end{align*}
	where the hat means that term is omitted. Set $f_i:=e_i\odot\dots\odot\hat{e_i}\odot\dots\odot  e_n$ for each $i\in\{1,..., n\}$. Hence $L=[e_i)\oplus[f_i)$  for each $i\in\{1,..., n\}$. By Propositions \ref{11} and \ref{1222}, $[f_i)$ is a maximal filter  for each $i\in\{1,..., n\}$ and $f_i=e^*_i$. We show that $Max(L)=\{[f_1),..., [f_n)\}$. By above argument, $\{[f_1),..., [f_n)\}\subseteq Max(L)$. Now let $N\in Max(L)$. For each $i\in\{1,..., n\}$, we have $[e_i)\cap[f_i)=\{1\}\subseteq N$. Hence for each $i\in\{1,..., n\}$, we have either  $[e_i)\subseteq N$ or $[f_i)\subseteq N$. If there exists $i\in\{1,..., n\}$ such that $[f_i)\subseteq N$, then $[f_i)= N$ since $[f_i)$ is a maximal filter of $L$, and so $N\in \{[f_1),..., [f_n)\}$. Now assume that for all $i\in\{1,..., n\}$,   $[e_i)\subseteq N$, then $L= [e_1)\vee\dots\vee [e_n)\subseteq N$. Thus $L=N$, a contradiction. Therefore, $Max(L)=\{[f_1),..., [f_n)\}$ and so $L$ is semi-local.
	
	$8\Rightarrow 1)$ Assumption assume that $Max(L)=\{M_1,..., M_n\}$. Since $L$ is hyperarchimedean, we have $ Rad(L)=\{1\}$ by Theorem \ref{hyp1}. Now if $n=1$, then $Max(L)=\{M_1\}$ and $M_1= Rad(L)=\{1\}$, and so $L$ is simple and we are done. Now assume that $n\geq 2$. For each $i\in\{1,..., n\} $ set $N_i:=M_i\cap\dots\cap\hat{M_i}\cap\dots\cap M_n$, where the hat means that term is omitted. Hence for each $i\in\{1,..., n\} $, we have  $M_i\cap N_i= Rad(L)=\{1\}$ and $M_i\vee N_i=M_i\vee (M_i\cap\dots\cap\hat{M_i}\cap\dots\cap M_n)=\bigcap_{j\not= i}(M_i\vee M_j)=\bigcap_{j\not= i}L=L$. Hence $L=M_i\oplus N_i$ for each $i$. Thus by Proposition \ref{1222}, $N_i$ is a simple filter of $L$ for each $i\in\{1,..., n\}$. So for each $i\in\{1,..., n\} $ there exists $e_i\in  B (L)$ such that $N_i=[e_i)$ and $M_i=[e^*_i)$. Since $e^*_1\vee\dots\vee e^*_n\in \bigcap_{i=1}^nM_i= Rad(L)=\{1\}$, we have $e^*_i\vee\dots\vee e^*_n=1 $. Set $f:=e_i\odot\dots\odot e_n$, hence we have:
	\begin{align*}
	f=f\odot 1=&f\odot(e^*_1\vee\dots\vee e^*_n)\\=&(f\odot e^*_1)\vee\dots\vee (f\odot e^*_n)\\\leq & (e_1\odot e^*_1)\vee\dots\vee (e_n\odot e^*_n)=0.
	\end{align*}
	Hence $f=e_i\odot\dots\odot e_n=0$. Thus $L=N_1\vee\dots\vee N_n$ and we are done.
	
	$8\Rightarrow 9)$ It follows from Theorem \ref{hyp1}.
	
	$9\Rightarrow 8)$ Assume that $Max(L)=\{M_1,..., M_n\}$. Hence $ Rad(L)=M_1\cap\dots\cap M_n=\{1\}$. Now if $P\in   Spec(L)$, then $M_1\cap\dots\cap M_n=\{1\}\subseteq P$. Hence, there exists $i\in\{1,..., n\}$ such that $M_i\subseteq P$. Since $M_i\in Max(L)$, we have $M_i=P$. Thus, $P\in Max(L)$, and so $  Spec(L)=Max(L)$. Therefore, $L$ is hyperarchimedean by Theorem \ref{hyp}.
	
	$8\Rightarrow 10)$ Assume that $L$ is hyperarchimedean and semi-local. By Theorem \ref{hyp}, $  Spec(L)=Max(L)$, and so by Proposition \ref{2.700}(2), every proper filter of $L$ can be written  as a finite intersection of distinct maximal filters. Now let $P_1\cap\dots\cap P_n$ and $Q_1\cap\dots\cap Q_m$ be two finite intersections of distinct maximal filters that are equal to a proper filter $F$ of $L$. Then for each $i\in\{1,..., n\}$  since $Q_1\cap\dots\cap Q_m\subseteq P_i$, there exists a unique  $j_i\in\{1,..., m\}$ such that $Q_{j_i}\subseteq P_i$. Now since $P_i, Q_{j_i}\in Max(L)$, we have $Q_{j_i}= P_i$. Similarly,  for each $j\in\{1,..., m\}$  since $P_1\cap\dots\cap P_n\subseteq Q_j$, there exists a unique  $i_j\in\{1,..., n\}$ such that $P_{i_j}\subseteq Q_j$. Now since $P_{i_j}, Q_j\in Max(L)$, we have $Q_j= P_{i_j}$. Hence, there exists a one to one correspondence between two sets $\{P_1,..., P_n\}$ and $\{ Q_1,..., Q_m\}$ and correspondence elements are equal. Hence,  every proper filter of $L$ can be written uniquely, up to order, as a finite intersection of distinct maximal filters. 
	
	$10\Rightarrow 11)$ It is clear.  
	
	$11\Rightarrow 8)$ Assume that $\{1\}=P_1\cap\dots\cap P_n$  be a finite intersection of distinct maximal filters for the filter  $\{1\}$ of $L$. Now if $P\in  Spec(L)$, then from $P_1\cap\dots\cap P_n=\{1\}\subseteq P$ there exists  $i\in\{1,..., n\}$ such that $P_i\subseteq P$. Now since $P_i\in Max(L)$, we have  $P_i= P$. Hence $Max(L)\subseteq   Spec(L)\subseteq \{P_1,..., P_n\}\subseteq Max(L)$. This shows that $Max(L)=   Spec(L)= \{P_1,..., P_n\}$. Hence $L$ is hyperarchimedean and semi-local by Theorem \ref{hyp1}.
	
	$1\Leftrightarrow 12)$ It is clear by Corollary \ref{4.3}.
\end{proof}
In the following result, we show that in the finite case semi-simple residuated lattices  and hyperarchimedean residuated lattices are the same.
\begin{corollary}\label{4.5}
	Let $L$ be a semi-local (e.g., finite residuated lattices) residuated lattice. Then the following are equivalent:
	\begin{enumerate}
		\item $L$ is semi-simple;
		\item $L$ is hyperarchimedean;
		\item $ Rad(L)=\{1\}$.
		
	\end{enumerate}
\end{corollary}
\begin{proof}
	The corollary follows from Theorem \ref{4.4}.
\end{proof}
In the following example, we show that the condition \lq\lq$L$ is semi-local" in Corollary \ref{4.5} is necessary.
\begin{example}
	Let $\{0, 1\}$ be the two-point Boolean algebra. Set $A_i:=\{0, 1\}$ for each $i\in\mathbb{N}$ and consider the infinite
	direct products of residuated lattices $A_i's$, that is, the residuated lattice $L=\prod_{i\in\mathbb{N}}A_i$. The clearly $ Rad(L)=\{(1, 1,...)\}=\{1_L\}$, but $L$ is not semi-local. Hence $L$ is not semi-simple by Theorem \ref{4.4}. Thus Corollary \ref{4.5} is not true in general without the hypothesis that $L$ is semi-local.
\end{example}

We end this section with a comparison between a ring's semi-simplicity and the semi-simplicity of the residuated lattice of its ideals.
 
In the literature, there are some equivalent definitions for semi-simple rings, here we recall one of them. Recall that a commutative unitary ring $R$ is called \textit{semi-simple}, whenever it is a direct product of finitely many fields, see [17, Theorem 3.5.19]. 
 
 For a commutative unitary ring $R$, it is known that the
lattice of ideals $(Id(R),\cap ,+,\odot, \rightarrow ,0=\{0\},1=R)$ is a
residuated lattice and the order relation is $\subseteq $, see   [14]. Actually, for every $I,J\in Id\left( R\right)$:

\[I+J:=\{i+j,i\in I,j\in J\}\text{, }I\odot J:=\{\underset{k=1}{\overset{n}{\sum }}i_{k}j_{k},\text{ }i_{k}\in
I,j_{k}\in J\}\text{, and }I\rightarrow J:=\{x\in A,x\cdot I\subseteq J\}.\]
\begin{theorem}\label{ring}
Let $R$ be commutative unitary ring. Then the residuated lattice $(Id(R),\cap ,+,\odot, \rightarrow ,0=\{0\},1=R)$ is semi-simple.
\end{theorem}
\begin{proof}
Let $R:=R_1\times\dots\times R_n$ be a semi-simple ring, where $R_i's$
 are fields. It is well-knowen that $Id(R)=\{I_1\times\dots\times I_n\mid I_i\in Id(R_i)\text{  for each } i=1,..., n\}$. Now since $R_i$ is a field, we have $Id(R_i)=\{\{0_{R_i}\}, R_i\}$. Hence, $(Id(R),\cap ,+,\odot, \rightarrow ,0=\{0\},1=R)$ is a finite residuated lattice. Clearly, for each $I_1\times\dots\times I_n\in R$ and $m\in \mathbb{N}$, we have $(I_1\times\dots\times I_n)^m=I_1\times\dots\times I_n$, and $(I_1\times\dots\times I_n)^*=J_1\times\dots\times J_n$, where for each $i$
 \begin{equation*}
J_i = \left\{
\begin{array}{rl}
\{0_{R_i}\} & \text{if } I_i=R_i,\\
R_i & \text{if } I_i=\{0_{R_i}\}.
\end{array} \right.
\end{equation*} 
 
 Hence for $I_1\times\dots\times I_n\in Id(R)$, if there is $k\in \mathbb{N}$ such that $((I_1\times\dots\times I_n)^*)^k=0_{Id(R)}=\{0_R\}=\{(0_{R_1},..., 0_{R_n})\}$, then $I_1\times\dots\times I_n=R_1\times\dots\times R_n=R=1_{Id(R)}$. Thus, $ Rad(Id(R))=\{1_{Id(R)}\}$ by Theorem \ref{rad}. Therefore, the result follows from Corollary \ref{4.5}.
 \end{proof}
 The converse of Theorem \ref{ring} is not true in general.
 \begin{example}
 Let $K$ be a field, and consider the ring $R:=K[[x]]/\langle x^2\rangle$, the quotient of  the power series ring $K[[x]]$ over the field $K$ by the ideal generated by $x^2$. An easy ring theory argument shows that $Id(R)=\{\{0_R\}:=\langle x^2\rangle/\langle x^2\rangle, a:=\langle x\rangle/\langle x^2\rangle, 1:=R\}$  is a residuated lattice with Hasse diagram depicted in the following figure, and  the operations $\odot$ and $\rightarrow$ via the accompanying tables:		
	\begin{center}  
		\begin{minipage}{0.4\textwidth}
			\centering
			\[
			\xymatrix@=8pt{
				& & 1 \\
				&   &a \ar@{-}[u] \ar@{-}[d]&    \\ 
				& & 0 & }		\]
		\end{minipage}
		\begin{minipage}{0.58\textwidth}
			\centering
			\begin{tabular}{r|l l l }
				$\odot$ & 0 & $a$   &  1 \\\hline
				0& $0$ &$0$& $0$\\
				$a$ &$0$ &$0$ & $a$\\
				1 & $0$& $a$&$1$\\
			\end{tabular} 
			\hspace{1cm}
			\begin{tabular}{r|l l l }
				$\rightarrow$ & 0 & $a$ &  1 \\\hline
				0& $1$ &$1$& $1$\\
				$a$ & $a$ & $1$ & $1$\\
				1 & $0$& $a$ & $1$\\
			\end{tabular}
			
			\vspace{0.5cm}
		\end{minipage}
	\end{center}
 Now since $R$ is a local ring that is not a field, it is not a semi-simple ring, but since $(a^2)^*=1$, we have for each $n\in \mathbb{N}$, $((a^2)^*)^n=1^n=1\not= 0$, and hence $ Rad(Id(R))=\{1\}$. Therefore, by Corollary \ref{4.5}, $Id(R)$ is a semi-simple residuated lattice.
\end{example}



\begin{thebibliography}{99}

\bibitem{77} Ahadpanah, A., \& Torkzadeh, L. (2015). Fuzzy Noetherian and Artinian residuated lattice. \textit{University Politehnica of Bucharest Scientific Bulletin-Series A-Applied Mathematics and Physics, 77}(4), 23–32.

\bibitem{anderson} Anderson, F. W., \& Fuller, K. R. (1992). \textit{Rings and categories of modules}. Springer-Verlag.

\bibitem{5} Ebrahimi Atani, S., Dolati Pish Hesari, S., \& Khoramdel, M. (2018). Semisimple semirings with respect to co-ideals theory. \textit{Asian Journal of Mathematics, 11}(3), [13 pages].

\bibitem{eng} Engelking, R. (1989). \textit{General topology}. Heldermann Verlag.

\bibitem{nga1} Galatos, N., Jipsen, P., Kowalski, T., \& Ono, H. (2007). \textit{Residuated lattices: An algebraic glimpse at substructural logics}. Elsevier.

\bibitem{ggg0} Georgescu, G., Cheptea, D., \& Muresan, C. (2015). Algebraic and topological results on lifting properties in residuated lattices. \textit{Fuzzy Sets and Systems, 271}, 102–132.

\bibitem{pha} Hájek, P. (1998). \textit{Metamathematics of fuzzy logic}. Kluwer Academic Publishers.

\bibitem{10} Hebisch, U., \& Weinert, H. J. (2002). Semisimple classes of semirings. \textit{Algebra Colloquium, 9}(2), 177–196.

\bibitem{dem} Holdon, L. C. (2018). On ideals in De Morgan residuated lattices. \textit{Kybernetika, 54}(3), 443–475.

\bibitem{11} Katsov, Y., Nam, T. G., \& Tuyen, N. X. (2009). On subtractive semisimple semirings. \textit{Algebra Colloquium, 16}(3), 415–426.

\bibitem{lam} Lam, T. Y. (2001). \textit{A first course in noncommutative rings} (2nd ed.). Springer-Verlag.

\bibitem{66} Motamed, S., \& Moghaderi, J. (2012). Noetherian and artinian BL-algebras. \textit{Soft Computing, 18}, 419–429. https://doi.org/10.1007/s00500-012-0926-1

\bibitem{dpi2} Piciu, D. (2007). \textit{Algebras of fuzzy logic}. Editura Universitaria Craiova.

\bibitem{ttt} Tchoffo Foka, S. V., \& Tonga, M. (2022). Rings and residuated lattices whose fuzzy ideals form a Boolean algebra. \textit{Soft Computing, 26}, 535–539. https://doi.org/10.1007/s00500-021-06430-9

\bibitem{etu} Turunen, E. (1999). \textit{Mathematics behind fuzzy logic}. Physica-Verlag.

\bibitem{wa} Wang, F., \& Kim, H. (2016). \textit{Foundations of commutative rings and their modules}. Springer.

\bibitem{111} Ward, M., \& Dilworth, R. P. (1939). Residuated lattices. \textit{Transactions of the American Mathematical Society, 45}(3), 335–354.

\bibitem{14} Wisbauer, R. (1991). \textit{Foundations of module and ring theory}. Gordon and Breach.

\end{thebibliography}
\end{document}